\documentclass[11pt]{article}
\usepackage{amsmath, upgreek}
\usepackage{amssymb}
\usepackage{amscd}
\usepackage{xspace}
\usepackage{verbatim}
\usepackage{color}
\newcommand{\mysection}[1]{
\section{#1}\setcounter{equation}{0}}
\title{\bf Trace and boundary singularities of positive solutions of a class of quasilinear equations}
\author{{\bf Marie-Fran\c{c}oise Bidaut-V\'eron}\\
 {\bf Laurent V\'eron}\\[2mm]
}

\date{}
\begin{document}
 \maketitle


\newcommand{\txt}[1]{\;\text{ #1 }\;}
\newcommand{\tbf}{\textbf}
\newcommand{\tit}{\textit}
\newcommand{\tsc}{\textsc}
\newcommand{\trm}{\textrm}
\newcommand{\mbf}{\mathbf}
\newcommand{\mrm}{\mathrm}
\newcommand{\bsym}{\boldsymbol}
\newcommand{\scs}{\scriptstyle}
\newcommand{\sss}{\scriptscriptstyle}
\newcommand{\txts}{\textstyle}
\newcommand{\dsps}{\displaystyle}
\newcommand{\fnz}{\footnotesize}
\newcommand{\scz}{\scriptsize}
\newcommand{\be}{\begin{equation}}
\newcommand{\bel}[1]{\begin{equation}\label{#1}}
\newcommand{\ee}{\end{equation}}
\newcommand{\eqnl}[2]{\begin{equation}\label{#1}{#2}\end{equation}}
\newcommand{\barr}{\begin{eqnarray}}
\newcommand{\earr}{\end{eqnarray}}
\newcommand{\bars}{\begin{eqnarray*}}
\newcommand{\ears}{\end{eqnarray*}}
\newcommand{\nnu}{\nonumber \\}
\newtheorem{subn}{\name}
\renewcommand{\thesubn}{}
\newcommand{\bsn}[1]{\def\name{#1}\begin{subn}}
\newcommand{\esn}{\end{subn}}
\newtheorem{sub}{\name}[section]
\newcommand{\dn}[1]{\def\name{#1}}   
\newcommand{\bs}{\begin{sub}}
\newcommand{\es}{\end{sub}}
\newcommand{\bsl}[1]{\begin{sub}\label{#1}}
\newcommand{\bth}[1]{\def\name{Theorem}
\begin{sub}\label{t:#1}}
\newcommand{\blemma}[1]{\def\name{Lemma}
\begin{sub}\label{l:#1}}
\newcommand{\bcor}[1]{\def\name{Corollary}
\begin{sub}\label{c:#1}}
\newcommand{\bdef}[1]{\def\name{Definition}
\begin{sub}\label{d:#1}}
\newcommand{\bprop}[1]{\def\name{Proposition}
\begin{sub}\label{p:#1}}
\newcommand{\R}{\eqref}
\newcommand{\rth}[1]{Theorem~\ref{t:#1}}
\newcommand{\rlemma}[1]{Lemma~\ref{l:#1}}
\newcommand{\rcor}[1]{Corollary~\ref{c:#1}}
\newcommand{\rdef}[1]{Definition~\ref{d:#1}}
\newcommand{\rprop}[1]{Proposition~\ref{p:#1}}
\newcommand{\BA}{\begin{array}}
\newcommand{\EA}{\end{array}}
\newcommand{\BAN}{\renewcommand{\arraystretch}{1.2}
\setlength{\arraycolsep}{2pt}\begin{array}}
\newcommand{\BAV}[2]{\renewcommand{\arraystretch}{#1}
\setlength{\arraycolsep}{#2}\begin{array}}
\newcommand{\BSA}{\begin{subarray}}
\newcommand{\ESA}{\end{subarray}}
\newcommand{\BAL}{\begin{aligned}}
\newcommand{\EAL}{\end{aligned}}
\newcommand{\BALG}{\begin{alignat}}
\newcommand{\EALG}{\end{alignat}}
\newcommand{\BALGN}{\begin{alignat*}}
\newcommand{\EALGN}{\end{alignat*}}
\newcommand{\note}[1]{\textit{#1.}\hspace{2mm}}
\newcommand{\Proof}{\note{Proof}}
\newcommand{\qeda}{\hspace{10mm}\hfill $\square$}
\newcommand{\qed}{\\
${}$ \hfill $\square$}
\newcommand{\Remark}{\note{Remark}}
\newcommand{\modin}{$\,$\\[-4mm] \indent}
\newcommand{\forevery}{\quad \forall}
\newcommand{\set}[1]{\{#1\}}
\newcommand{\setdef}[2]{\{\,#1:\,#2\,\}}
\newcommand{\setm}[2]{\{\,#1\mid #2\,\}}
\newcommand{\mt}{\mapsto}
\newcommand{\lra}{\longrightarrow}
\newcommand{\lla}{\longleftarrow}
\newcommand{\llra}{\longleftrightarrow}
\newcommand{\Lra}{\Longrightarrow}
\newcommand{\Lla}{\Longleftarrow}
\newcommand{\Llra}{\Longleftrightarrow}
\newcommand{\warrow}{\rightharpoonup}
\newcommand{
\paran}[1]{\left (#1 \right )}
\newcommand{\sqbr}[1]{\left [#1 \right ]}
\newcommand{\curlybr}[1]{\left \{#1 \right \}}
\newcommand{\abs}[1]{\left |#1\right |}
\newcommand{\norm}[1]{\left \|#1\right \|}
\newcommand{
\paranb}[1]{\big (#1 \big )}
\newcommand{\lsqbrb}[1]{\big [#1 \big ]}
\newcommand{\lcurlybrb}[1]{\big \{#1 \big \}}
\newcommand{\absb}[1]{\big |#1\big |}
\newcommand{\normb}[1]{\big \|#1\big \|}
\newcommand{
\paranB}[1]{\Big (#1 \Big )}
\newcommand{\absB}[1]{\Big |#1\Big |}
\newcommand{\normB}[1]{\Big \|#1\Big \|}
\newcommand{\produal}[1]{\langle #1 \rangle}

\newcommand{\thkl}{\rule[-.5mm]{.3mm}{3mm}}
\newcommand{\thknorm}[1]{\thkl #1 \thkl\,}
\newcommand{\trinorm}[1]{|\!|\!| #1 |\!|\!|\,}
\newcommand{\bang}[1]{\langle #1 \rangle}
\def\angb<#1>{\langle #1 \rangle}
\newcommand{\vstrut}[1]{\rule{0mm}{#1}}
\newcommand{\rec}[1]{\frac{1}{#1}}
\newcommand{\opname}[1]{\mbox{\rm #1}\,}
\newcommand{\supp}{\opname{supp}}
\newcommand{\dist}{\opname{dist}}
\newcommand{\myfrac}[2]{{\displaystyle \frac{#1}{#2} }}
\newcommand{\myint}[2]{{\displaystyle \int_{#1}^{#2}}}
\newcommand{\mysum}[2]{{\displaystyle \sum_{#1}^{#2}}}
\newcommand {\dint}{{\displaystyle \myint\!\!\myint}}
\newcommand{\q}{\quad}
\newcommand{\qq}{\qquad}
\newcommand{\hsp}[1]{\hspace{#1mm}}
\newcommand{\vsp}[1]{\vspace{#1mm}}
\newcommand{\ity}{\infty}
\newcommand{\prt}{\partial}
\newcommand{\sms}{\setminus}
\newcommand{\ems}{\emptyset}
\newcommand{\ti}{\times}
\newcommand{\pr}{^\prime}
\newcommand{\ppr}{^{\prime\prime}}
\newcommand{\tl}{\tilde}
\newcommand{\sbs}{\subset}
\newcommand{\sbeq}{\subseteq}
\newcommand{\nind}{\noindent}
\newcommand{\ind}{\indent}
\newcommand{\ovl}{\overline}
\newcommand{\unl}{\underline}
\newcommand{\nin}{\not\in}
\newcommand{\pfrac}[2]{\genfrac{(}{)}{}{}{#1}{#2}}

\def\ga{\alpha}     \def\gb{\beta}       \def\gg{\gamma}
\def\gc{\chi}       \def\gd{\delta}      \def\ge{\epsilon}
\def\gth{\theta}                         \def\vge{\varepsilon}
\def\gf{\phi}       \def\vgf{\varphi}    \def\gh{\eta}
\def\gi{\iota}      \def\gk{\kappa}      \def\gl{\lambda}
\def\gm{\mu}        \def\gn{\nu}         \def\gp{\pi}
\def\vgp{\varpi}    \def\gr{\rho}        \def\vgr{\varrho}
\def \gs{\sigma}     \def\vgs{\varsigma}  \def\gt{\tau}
\def\gu{\upsilon}   \def\gv{\vartheta}   \def\gw{\omega}
\def\gx{\xi}        \def\gy{\psi}        \def\gz{\zeta}
\def\Gg{\Gamma}     \def\Gd{\Delta}      \def\Gf{\Phi}
\def\Gth{\Theta}
\def\Gl{\Lambda}    \def\Gs{\Sigma}      \def\Gp{\Pi}
\def\Gw{\Omega}     \def\Gx{\Xi}         \def\Gy{\Psi}

\def\CS{{\mathcal S}}   \def\CM{{\mathcal M}}   \def\CN{{\mathcal N}}
\def\CR{{\mathcal R}}   \def\CO{{\mathcal O}}   \def\CP{{\mathcal P}}
\def\CA{{\mathcal A}}   \def\CB{{\mathcal B}}   \def\CC{{\mathcal C}}
\def\CD{{\mathcal D}}   \def\CE{{\mathcal E}}   \def\CF{{\mathcal F}}
\def\CG{{\mathcal G}}   \def\CH{{\mathcal H}}   \def\CI{{\mathcal I}}
\def\CJ{{\mathcal J}}   \def\CK{{\mathcal K}}   \def\CL{{\mathcal L}}
\def\CT{{\mathcal T}}   \def\CU{{\mathcal U}}   \def\CV{{\mathcal V}}
\def\CZ{{\mathcal Z}}   \def\CX{{\mathcal X}}   \def\CY{{\mathcal Y}}
\def\CW{{\mathcal W}} \def\CQ{{\mathcal Q}}
\def\BBA {\mathbb A}   \def\BBb {\mathbb B}    \def\BBC {\mathbb C}
\def\BBD {\mathbb D}   \def\BBE {\mathbb E}    \def\BBF {\mathbb F}
\def\BBG {\mathbb G}   \def\BBH {\mathbb H}    \def\BBI {\mathbb I}
\def\BBJ {\mathbb J}   \def\BBK {\mathbb K}    \def\BBL {\mathbb L}
\def\BBM {\mathbb M}   \def\BBN {\mathbb N}    \def\BBO {\mathbb O}
\def\BBP {\mathbb P}   \def\BBR {\mathbb R}    \def\BBS {\mathbb S}
\def\BBT {\mathbb T}   \def\BBU {\mathbb U}    \def\BBV {\mathbb V}
\def\BBW {\mathbb W}   \def\BBX {\mathbb X}    \def\BBY {\mathbb Y}
\def\BBZ {\mathbb Z}   \def\BBQ {\mathbb Q}

\def\GTA {\mathfrak A}   \def\GTB {\mathfrak B}    \def\GTC {\mathfrak C}
\def\GTD {\mathfrak D}   \def\GTE {\mathfrak E}    \def\GTF {\mathfrak F}
\def\GTG {\mathfrak G}   \def\GTH {\mathfrak H}    \def\GTI {\mathfrak I}
\def\GTJ {\mathfrak J}   \def\GTK {\mathfrak K}    \def\GTL {\mathfrak L}
\def\GTM {\mathfrak M}   \def\GTN {\mathfrak N}    \def\GTO {\mathfrak O}
\def\GTP {\mathfrak P}   \def\GTR {\mathfrak R}    \def\GTS {\mathfrak S}
\def\GTT {\mathfrak T}   \def\GTU {\mathfrak U}    \def\GTV {\mathfrak V}
\def\GTW {\mathfrak W}   \def\GTX {\mathfrak X}    \def\GTY {\mathfrak Y}
\def\GTZ {\mathfrak Z}   \def\GTQ {\mathfrak Q}

\font\Sym= msam10 
\def\SYM#1{\hbox{\Sym #1}}
\newcommand{\bdw}{\prt\Gw\xspace}
\medskip

\nind {\it 2010 Mathematics Subject Classification:} 35J62-35J66-35J75-31C15\smallskip

\nind{\it  Keywords:} Elliptic equations, boundary singularities, Bessel capacities, measures, supersolutions, subsolutions.
\medskip

\noindent {\it \small A Juan-Luis por su 75 cumplea\~nos. Cuarenta y seis a\~nos de amistad, respeto y admiraci\'{o}n}

{\abstract We study positive functions satisfying (E)$\;-\Gd u+m|\nabla u|^q-u^p=0$ in a domain $\Gw$ or in $\BBR^{_N}_+$ when $p>1$ and $1<q<2$. We give sufficient conditions for the existence of a solution to (E) with a nonnegative measure $\gm$ as boundary data; these conditions are expressed in terms of Bessel capacities on the boundary. We also study removability of boundary singular sets, and solutions with an isolated singularity on $\prt\Gw$. The different results depend on two critical exponents for $p=p_c:=\frac{N+1}{N-1}$ and for $q=q_c:=\frac{N+1}{N}$, and on the sign of $q-\frac{2p}{p+1}$.
\tableofcontents
\date{}
\maketitle

\medskip

\mysection{Introduction}
In this article we study the boundary behaviour of positive solutions of the following class of quasilinear elliptic equations 
\begin{equation}\label{An1}
\BA{lll}
-\Gd u+m|\nabla u|^{q}-|u|^{p-1}u=0
\EA
\end{equation}
in  a domain $G$ of $\BBR^N$ which can be either $\BBR^N$, or $\BBR^N\setminus\{0\}$, or $\BBR^N_+$, or a bounded domain $\Gw$ with smooth boundary $\prt\Gw$, according the type of phenomenon we are interested in. We assume that $p,q>1$ and $m\geq 0$. We also consider the associated measure boundary data problem 
\begin{equation}\label{An2}
\BA{lll}
-\Gd u+m|\nabla u|^{q}-|u|^{p-1}u=0\qquad&\text{in }\; \Gw\\
\phantom{-\Gd u+m|\nabla u|^{q}-|u|^{p-1}}
u=\gm\qquad&\text{in }\; \prt\Gw,
\EA
\end{equation}
in the case where $G=\Gw$ and $\gm$ is a positive Radon measure on $\prt\Gw$.\\

The wide variety of phenomena that exhibit the solutions of equation $(\ref{An1})$ comes from the opposition between the forcing term $|u|^{p-1}u$ and the reaction term $m|\nabla u|^{q}$. Furthermore,  in the specific case $q=\frac{2p}{p+1}$, the order of magnitude of the forcing and the reaction is the same, therefore the value of the coefficient $m$ plays  a fundamental role. This is due to the equivariance of the equation $(\ref{An1})$ under the transformation $u\mapsto T_\ell[u]$ defined by 
 \begin{equation}\label{An3}
T_\ell[u](x)=\ell^{\frac{2}{p-1}}u(\ell x)\quad\text{where }\ell>0. 
\end{equation}
  This equation has been introduced by Chipot and Weissler in \cite{ChWe} in a parabolic setting. They also studied the one dimensional case of $(\ref{An1})$.  Later on Serrin and Zou published two deep articles \cite{SeZo-1},  \cite{SeZo-2}
   where they concentrate on the existence of radial ground states, introducing unexpected energy functions. In \cite{SeYaZo} they conduct a series of numerical experimentations showing the extreme complexity of this equation, even in the radial case, and many deep questions that they raised are still  unanswered. More recently, Alarc\'on, Garc\'ia-Meli\'an and Quaas proved several non-existence results of supersolutions in an exterior domain of a large class of equations containing in particular $(\ref{An1})$. Their results pointed out the role of some critical exponents, $p=\frac{N}{N-2}$, $p=\frac{N+2}{N-2}$ and $q=\frac{N}{N-1}$ as well as $q=\frac{2p}{p+1}$. A priori estimates of solutions have been obtained  in \cite{PQS} in the case $q<\frac{2p}{p+1}$ and $p<\frac{N+2}{N-2}$, and then extended in \cite{BVGHV1} to the case $q=\frac{2p}{p+1}$ and $p<\frac{N+2}{N-2}$ under a condition of smallness of $m$ by a completely different method. 
The regular Dirichlet problem has been investigated in \cite{Ruiz} in the subcritical case $p<\frac N{N-2}$ and $q<\frac {2p}{p+1}$, and even extended to the $m$-Laplace equation, always in the corresponding subcritical case, but to our knowledge, nothing has already been published concerning the boundary behaviour of singular solutions and the associated Dirichlet problem with measure data. The aim of this article is to fulfill some gaps in the knowledge of the properties of this equation, emphasizing the connection with an acurate description of the boundary behaviour. \smallskip

We first prove an a priori estimate for positive solutions of 
 \begin{equation}\label{An4}
\BA{lll}
-\Gd u+m|\nabla u|^{q}-u^p=0\qquad&\text{in }\Gw\\\phantom{-\Gd +m|\nabla u|^{q}-u^p}
u=0&\text{on }\prt \Gw\setminus\{0\}.
\EA
\end{equation}
when $q\leq \frac{2p}{p+1}$. We set
 \begin{equation}\label{Ann4}
\ga=\frac{2}{p-1}.
\end{equation}
\bth{BS1} Let $\Gw$ be a bounded smooth domain such that $0\in\prt\Gw$. Suppose $1<p<\frac{N+2}{N-2}$ and either $1<q<\frac{2p}{p+1}$ and $m>0$, or $q=\frac{2p}{p+1}$ and $0<m<\ge_0$ for some $\ge_0>0$ depending on $N$ and $p$. Then there exists a constant $c=c(N,p,\Gw)>0$ such that if $u$ is a positive solution 
of $(\ref{An4})$, it satisfies
 \begin{equation}\label{An5}
\BA{lll}
u(x)\leq c|x|^{-\ga}\qquad\text{for all }x\in\overline\Gw\setminus\{0\},
\EA\ee
and
 \begin{equation}\label{An6}
\BA{lll}
\myfrac{u(x)}{\gr(x)}+|\nabla u(x)|\leq c|x|^{-\ga-1}\quad\text{for all }x\in\overline\Gw\setminus\{0\},
\EA\ee
where $\gr(x)=\dist(x,\prt\Gw)$.
\es

Thanks to this estimate we can describe the behaviour of positive functions satisfying $(\ref{An4})$.  For this purpose we say that the bounded open set $\Gw\subset\BBR^N$ is in {\it normal position with respect to $0\in\prt\Gw$} if $\prt\Gw$ is tangent to $\prt\BBR_+^{N}$ at $x=0$ and if $x_{_N}>0$ is the normal inward direction to $\prt\Gw$. We set $\prt B_{1\,+}:=\BBR^N_+\cap\prt B_1$, identified with $S^{N-1}_+:=S^{N-1}\cap \BBR_+^{N}$ in spherical coordinates $(r,s)$. In the sequel we denote by $\Gd'$ the Laplace-Beltrami operator on $S^{N-1}$ and by $\nabla'$ the covariant gradient identified with the tangential gradient to $\prt B_1$.
\bcor{BSS1} Let $\Gw$ be a bounded smooth domain in normal position with respect to $0\in\prt\Gw$. Suppose $1<p<\frac{N+2}{N-2}$, $1<q<\frac{2p}{p+1}$ and $m>0$, and $u$ is a positive solution of $(\ref{An1})$ vanishing on $\prt\Gw\setminus \{0\}$, then either $u$ can be extended as a continuous function in $\overline\Gw$, or one of the following situations oocurs.\smallskip

\nind 1- If $1<p<\frac{N+1}{N-1}$, there exists $k> 0$ such that 
 \begin{equation}\label{An7}\displaystyle
\BA{lll}\displaystyle
\lim_{\tiny\BA{lll}x\in\Gw\\x\to 0\EA}\myfrac{u(x)}{\BBP_{\Gw}(x,0)}=k,
\EA\ee
where $\BBP_{\Gw}$ is the Poisson kernel in $\Gw$ with asymptotics given in $(\ref{An7*})$.\smallskip 

\nind 2- If $p=\frac{N+1}{N-1}$,
 \begin{equation}\label{An8}
\BA{lll}\displaystyle
\lim_{\tiny\BA{lll}x\in\Gw\\x\to 0\\
\frac{x}{|x|}\to s\EA}|x|^{N-1}\left(\ln\frac1{|x|}\right)^{\frac{N-1}{2}}u(x)=\gl_{_N}\phi_1( s),
\EA\ee
uniformly on any compact set of $S^{N-1}_+$, where  $\phi_1$ is the first eigenfunction of $-\Gd'$ in $W^{1,2}_0(S^{N-1}_+)$ with maximum $1$ (actually $\phi_1(x/|x|)=\sin(x_N/|x|)$), and $\gl_{_N}$ is a positive constant depending only on $N$.\smallskip 

\nind 3- If $\frac{N+1}{N-1}<p<\frac{N+2}{N-2}$ 
\begin{equation}\label{An9}
\BA{lll}\displaystyle
\lim_{\tiny\BA{lll}x\in\Gw\\x\to 0\\
\!\frac{x}{|x|}\to s\EA}|x|^{\ga}u(x)=\psi( s),
\EA\ee
with $\ga=\frac{2}{p-1}$, uniformly on any compact set of $S^{N-1}_+$, where  $\psi$ is the unique positive solution of 
\begin{equation}\label{An10}\BA{lll}
-\Gd'\psi+\ga(N-2-\ga)\psi-\psi^p=0\qquad&\text{in }\;S^{N-1}_+\\
\phantom{-\Gd'+\ga(N-2-\ga)\psi-\psi^p}
\psi=0\qquad&\text{on }\;\prt S^{N-1}_+.
\EA\ee
\es

A direct computation by matching asymptotic expansion shows that if $\Gw$ is in normal position at $0\in\prt\Gw$ the Poisson kernel has the following asymptotic expression near $x=0$
 \begin{equation}\label{An7*}\BBP_{\Gw}(x,0)=c_N|x|^{1-N}\left(\sin\left(\tfrac{x_N}{|x|}\right)+O\left(\tfrac{x^2_N}{|x|^2}\right)\right)
 \quad\text{as }x\in\Gw,\;x\to 0,
\ee
for some explicit constant $c_N$.\\
In case 1, a solution which satisfies $(\ref{An7})$ is actually a weak solution of 
 \begin{equation}\label{An11}
\BA{lll}
-\Gd u+m|\nabla u|^{q}-u^p=0\qquad&\text{in }\Gw\\\phantom{-\Gd +m|\nabla u|^{q}-u^p}
u=k\gd_0&\text{in }\CD'(\prt \Gw).
\EA
\end{equation}
where  $\gd_0$ is the Dirac measure at $0$. A solution which satisfies 
$(\ref{An8})$ or $(\ref{An9})$ is a weak solution of $(\ref{An1})$ in $\Gw$ with zero boundary value in the sense of distributions in $\prt\Gw$, and this property still holds even when $\frac{N+1}{N-1}\leq p<\frac{N+2}{N-2}$ and 
$q=\frac{2p}{p+1}$. \smallskip

The proof of \rcor{BSS1} is based upon the fact that if $1<q<\frac{2p}{p+1}$, the a priori estimates of \rth{BS1} imply that problem $(\ref{An4})$ is a perturbation of 
 \begin{equation}\label{An12}
\BA{lll}
-\Gd u-u^p=0\qquad&\text{in }\Gw\\\phantom{-\Gd -u^p}
u=0&\text{on }\prt \Gw\setminus\{0\}
\EA
\end{equation}
near $x=0$, a problem which has been thoroughly studied in \cite{BVPV}. When $q=\frac{2p}{p+1}$, it is a consequence of the invariance of $(\ref{An1})$
 under the transformations $T_\ell$ that there could exist invariant solutions $u$ which are the ones such that $T_\ell[u]=u$ for any $\ell>0$. We first consider self-similar solutions in whole $\BBR^N$. Using spherical coordinates $(r,s)\in \BBR_+\ti S^{N-1}$, these self-similar solutions have the form 
   \begin{equation}\label{A5-1}
\BA{lll}
u(x)=u(r,s)=r^{-\ga}\gw(s),
\EA
\end{equation}
where $\ga$ is defined in $(\ref{Ann4})$. Then $\gw$ satisfies
   \begin{equation}\label{A6}
\BA{lll}
-\Gd'\gw+\ga(N-2-\ga)\gw+m\left(\ga^2\gw^2+|\nabla'\gw|^2\right)^{\frac{p}{p+1}}-|\gw|^{p-1}\gw=0\EA
\end{equation}
in $S^{N-1}$. Constant  solutions are roots of the function
   \begin{equation}\label{A7}
\BA{lll}
\CP_m(X)=\ga(N-2-\ga)X+m\ga^{\frac{2p}{p+1}}X^{\frac{2p}{p+1}}-|X|^{p-1}X.
\EA
\end{equation}
In the study of the variations of $\CP_m$ on $\BBR$ the following constant, defined if $p<\frac{N}{N-2}$, plays an important role
   \begin{equation}\label{A8}
\BA{lll}
m^*=(p+1)\left(\myfrac{N-p(N-2)}{2p}\right)^{\frac p{p+1}}.
\EA
\end{equation}
Concerning the self-similar solutions in $\BBR^N$ we recall the result stated without proof in \cite[Prop. 6.1]{BVGHV1},
\bprop{rad} Assume $N\geq 2$. \smallskip

\nind (i) If $N\geq 3$, $m> 0$  and $p\geq \frac{N}{N-2}$ there exists a unique constant positive solution $X_{{m}}$ to $(\ref{A6})$.\smallskip

\nind (ii) If $N\geq 2$, $1<p<\frac{N}{N-2}$  and $m> m^*$ there exist two  constant positive solutions $0<X_{{1,m}}<X_{{2,m}}$ to $(\ref{A6})$.\smallskip

\nind (iii) If $N\geq 2$ $1<p<\frac{N}{N-2}$ and $m= m^*$ there exists a unique constant positive solution $X_{m*}$ to $(\ref{A6})$.\smallskip

\nind (iv) If $N\geq 2$, $1<p<\frac{N}{N-2}$ and $0<m< m^*$ there exists no constant positive solution to $(\ref{A6})$.
\es

A more complete study of equation $(\ref{A6})$ and its role in the description of isolated singularities is developed in the forthcomming paper \cite{BVGHV22}.\smallskip

When the domain $G$ in which we consider equation $(\ref{An1})$ has a non-empty boundary (in the sequel, either $G=\BBR^N_+$ or $G$ is a smooth bounded domain that we denote by $\Gw$), it is natural to study solutions of $(\ref{An1})$  with an isolated singularity lying on the boundary. The understanding of boundary singularities is conditioned by the knowledge of positive self-similar solutions in $\BBR^N_+$ vanishing on $\prt\BBR^N_+$ except at $x=0$. They are solutions of 
 \begin{equation}\label{A9}
\BA{lll}
-\Gd'\gw+\ga(N-2-\ga)\gw+m\left(\ga^2\gw^2+|\nabla'\gw|^2\right)^{\frac{p}{p+1}}-\gw^p=0\quad&\text{in }S^{N-1}_+\\\phantom{-\Gd'\gw+\ga(N-2-\ga)+m\left(\ga^2\gw^2+|\nabla'\gw|^2\right)^{\frac{p}{p+1}}-\gw^p}
\gw=0&\text{on }\prt S^{N-1}_+.
\EA
\end{equation}
There, the critical value for $p$ is $\frac{N+1}{N-1}$. The main result concerning problem $(\ref{A9})$ states as follows,

\bth{Exist} Let $1<p<\frac{N+1}{N-1}$.\\
\nind 1-For any $m\geq m^*$ there exists at least one positive solution $\gw_{m}$ to $(\ref{A9})$.\\
\nind 2- There exists $m_p \in (0,m^*)$ such that for any $0<m\leq m_p$ there exists no positive solution to $(\ref{A9})$.
\es
The value of $m_p$ is explicit. \smallskip

In the next section of this article we study problem $(\ref{An2})$. We denote by $L^1_\gr(\Gw)$ the space of measurable functions $u$ in $\Gw$ such that 
$u\gr\in L^1(\Gw)$ where $\gr(x)=\dist(x,\prt\Gw)$.
\bdef{weaksol} Let $p,q>0$, $m\in\BBR$ and $\gm$ be a bounded measure on $\prt\Gw$. We say that a nonnegative Borel function $u$ defined in $\Gw$ is a weak solution of  $(\ref{An2})$ if $u\in L^1(\Gw)$, $u^p\in L_\gr^1(\Gw)$, $|\nabla u|^{q}\in L^1_\gr(\Gw)$ and 
\begin{equation}\label{D1-1}
\BA{lll}
\myint{\Gw}{}\left(-u\Gd\gz+\left(m|\nabla u|^{q}-u^p\right)\gz\right) dx=-\myint{\prt\Gw}{}\myfrac{\prt\gz}{\prt{\bf n}}d\gm, 
\EA\ee
for all $\gz\in\BBX(\Gw):=\left\{\gz\in C^1_c(\overline\Gw):\Gd \gz\in L^{\infty}(\Gw)\right\}$.
\es

For $a>0$ and $1<b<\infty$, we denote by $Cap^{\prt\Gw}_{a,b}$ the Bessel capacity on $\prt\Gw$. It is defined by local 
charts (see e.g. \cite{MV01}). Our main existence result is the following.
\bth{meas} Let $p>1$, $1<q<2$ and $m>0$. Assume $\gm$ is a nonnegative measure on $\prt\Gw$. If $\gm$ satisfies
\begin{equation}\label{D6}
\BA{lll}
\gm(K)\leq C_3\min\left\{Cap^{\prt\Gw}_{\frac{2-q}{q},q'}(K),Cap^{\prt\Gw}_{\frac{2}{p},p'}(K), \right\}\quad\text{for any compact set } K\subset\prt\Gw,
\EA\ee
then one can find $\ge_3>0$ such that for any $0<\ge\leq \ge_3$ there exists a weak solution $u$ to problem $(\ref{An2})$ with $\gm$ replaced by $\ge\gm$.
\es
The main idea for proving this result is to associate to $(\ref{An1})$ the two problems
\begin{equation}\label{An2-1}
\BA{lll}
-\Gd v+m|\nabla v|^{q}=0\qquad&\text{in }\; \Gw\\
\phantom{-\Gd+m|\nabla v|^{q}}
v=\gm\qquad&\text{in }\; \prt\Gw,
\EA
\end{equation}
and
\begin{equation}\label{An2-2}
\BA{lll}
-\Gd w-w^p=0\qquad&\text{in }\; \Gw\\
\phantom{-\Gd -w^p}
w=\gm\qquad&\text{in }\; \prt\Gw.
\EA
\end{equation}
We show that when $(\ref{D8})$ holds these two problems admit positive solutions respectively $v_\gm$ and $w_\gm$,
 such that $0<v_\gm<w_\gm$ (with $\gm$ replaced by $\ge\gm$) which both satisfy the boundary trace relation as it is introduced in 
 \cite{MV00},
\begin{equation}\label{D01}
\BA{lll}\displaystyle
\lim_{\gd\to 0}\myint{\{\gr(x)=\gd\}}{}v_\gm ZdS=\lim_{\gd\to 0}\myint{\{\gr(x)=\gd\}}{}w_\gm ZdS=\myint{\prt\Gw}{}Z d\gm,
\EA\ee
for any $Z\in C(\overline\Gw)$. Since $v_\gm$ and $w_\gm$ are respectively a subsolution and a supersolution of $(\ref{An1})$, we derive the existence of a solution $u$ of $(\ref{An1})$ in $\Gw$ which satisfies also the boundary trace relation $(\ref{D01})$. This approach is linked to the dynamical construction of the boundary trace developed in \cite{MV02}. 
As an easy consequence of \rth{meas} we have the following result.

\bcor{cormeas1} Let $p>1$, $1<q<2$ and $m>0$. If $\gm$ is a nonnegative measure on $\prt\Gw$ there exists a 
positive weak solution to $(\ref{An2})$ with $\gm$ replaced by $\ge\gm$ under the following conditions.\smallskip

\nind 1- If $1<p<\frac{N+1}{N-1}$ and $1<q<\frac{N+1}{N}$. \smallskip

\nind 2- If $p\geq \frac{N+1}{N-1}$, $1<q<\frac{N+1}{N}$ and $\gm$ satisfies
 \begin{equation}\label{D9}
\BA{lll}
\gm(K)\leq C_2Cap^{\prt\Gw}_{\frac{2}{p},p'}(K)\quad\text{for any compact set } K\subset\prt\Gw.
\EA\ee

\nind 3- If $1<p<\frac{N+1}{N-1}$, $\frac{N+1}{N}\leq q<2$ and $\gm$ satisfies
 \begin{equation}\label{D9*}
\BA{lll}
\gm(K)\leq C_1Cap^{\prt\Gw}_{\frac{2-q}{q},q'}(K)\quad\text{for any compact set } K\subset\prt\Gw.
\EA\ee
\es

A more delicate corollary is based upon relations between Bessel capacities.

\bcor{cormeas2} Let $p>1$, $1<q<2$ and $m>0$. If $\gm$ is a nonnegative measure on $\prt\Gw$, there exists a 
positive weak solution to $(\ref{An2})$ with $\gm$ replaced by $\ge\gm$ under the following conditions.\smallskip

\nind 1- If  $\frac{N+1}{N}\leq q<\frac{2p}{p+1}$, when 
 \begin{equation}\label{D9-0}
\BA{lll}
\gm(K)\leq C_3Cap^{\prt\Gw}_{\frac{2}{p},p'}(K)\quad\text{for any compact set } K\subset\prt\Gw.
\EA\ee

\nind 2- If $p\geq\frac{N+1}{N-1}$ and  $q\geq\frac{2p}{p+1}$, when 
 \begin{equation}\label{D9-1}
\BA{lll}
\gm(K)\leq C_4Cap^{\prt\Gw}_{\frac{2-q}{q},q'}(K)\quad\text{for any compact set } K\subset\prt\Gw.
\EA\ee
\es

It is noticeable that the results of \rcor{cormeas1} and \rcor{cormeas2} cover the full range of exponents $(p,q)\in (1,\infty)\ti(1,2)$. The sufficient conditions of \rth{meas} are stronger than the necessary conditions which are obtained below.
\bth{meas-nec} Let $p>1$, $1<q<2$ and $m>0$. Assume there exists a nonnegative solution $u$ of problem $(\ref{An4})$ for some $\gm\in\mathfrak M_+(\prt\Gw)$. Then $\gm$ satisfies 
\begin{equation}\label{D6-1}
\BA{lll}
Cap^{\prt\Gw}_{\frac{2-q}{q},q'}(K)=0\Longrightarrow \gm(K)=0\quad\text{if $K\subset\prt\Gw$ is a compact set},
\EA\ee
and
\begin{equation}\label{D6-2}
\BA{lll}
Cap^{\prt\Gw}_{\frac{2}{p},p'}(K)=0\Longrightarrow \gm(K)=0\quad\text{if $K\subset\prt\Gw$ is a compact set}.
\EA\ee
\es

In the last section we study the boundary trace of positive solutions of $(\ref{An1})$. The notion of boundary trace is classical in harmonic analysis in the framework of bounded Borel measures. It has been extended to semilinear elliptic equations by Marcus and V\'eron in \cite{MV00}, \cite{MV01}, \cite{MV02} with general Borel measures as a natural framework for the boundary trace.
\bdef{bdrt} Let $\Gw\subset\BBR^N$ be a smooth bounded domain, $p>1$, $1<q<p$, $m>0$ and $\CO$ a relatively open subset of $\prt\Gw$. We say that a positive solution  $u$ of  $(\ref{An1})$ in $\Gw$ admits a boundary trace on $\CO$, denoted by $Tr_\CO(u)$, if there exist a relatively open subset $\CR(u)$ of $\CO$ and 
a nonnegative Radon measure $\gm$ on $\CR(u)$ such that 
\begin{equation}\label{D080}
\BA{lll}\displaystyle
\lim_{\gd\to 0}\myint{\{\gr(x)=\gd\}}{}u ZdS=\myint{\prt\Gw}{}Z d\gm
\EA\ee
for every $Z\in C(\overline\Gw)$ such that $supp (Z\lfloor_{\CO})\subset \CR(u)$, and if for every $z\in \CS(u):=\CO\setminus\CR(u)$ and any $\ge>0$, there holds
\begin{equation}\label{D09}
\BA{lll}\displaystyle
\lim_{\gd\to 0}\myint{\{\gr(x)=\gd\}\cap B_\ge(z)}{}udS=\infty.
\EA\ee
The boundary trace $Tr_\CO(u)$ is represented by the couple $(\CS(u),\gm)$ or equivalently by the outer Borel measure $\gm^*_\CO$ on $\CO$ defined as follows: 
\begin{equation}\label{D10}
\gm^*_\CO(\gz)=\left\{\BA{lll}\displaystyle\myint{\CR(u)}{}\gz d\gm&\;\forall\gz\in C^\infty(\prt\Gw)\text{ s.t. }\, supp(\gz)\subset\CR(u)\\[3mm]
\infty&\;\forall\gz\in C^\infty(\prt\Gw)\text{ s.t. }\, supp(\gz)\cap\CS(u)\neq\emptyset,\,\gz\geq 0. 
\EA\right.\ee
\es

It is easy to prove that if a compact set $K\subset \overline\Gw$ is such that $u^p+|\nabla u|^q\in L^1_{\gr}(K)$, then 
$K\cap\prt\Gw\subset\CR(u)$. 
We first give a result where the trace is always a nonnegative Radon measure. The meaning of this result is that the absorption term is dominated by the reaction term and the 
solution behaves like a superharmonic function.

\bth{Tra0}Let $\Gw$ be a bounded smooth domain and $p>1$. Assume  either $1<q<\frac {2p}{p+1}$ and $m>0$, or 
$q=\frac{2p}{p+1}$ and $0<m\leq m_1$ for some $m_1>0$ depending on $N$ and $p$. If $u$ is a positive solution of $(\ref{An1})$ in $\Gw$, then $u\in L^1(\Gw)$, $u^{p}+|\nabla u|^{q}\in L^1_\gr(\Gw)$ and there exists a nonnegative Radon measure $\gm$ on $\prt\Gw$ such that $u$ is a solution of  $(\ref{An2})$.
\es

The boundary trace of a positive solution of $(\ref{An1})$ may not be a Radon measure, for example, if 
$\frac{N+1}{N-1}<p<\frac{N}{N-2}$, $q=\frac{2p}{p+1}$, $m\geq m^*$ and $u$ is the restriction to $\Gw$ of a radial singular solution obtained in \rprop{rad}. In that case
\begin{equation}\label{D02}
\BA{lll}\displaystyle
\lim_{\gd\to 0}\myint{\{\gr(x)=\gd\}}{}u ZdS=\infty,
\EA\ee
for any $Z\in C_+(\overline\Gw)$, such that $Z(0)>0$. 
We have the following result.
\bth{obs1} Assume $p>1$, $1<q<p$ and $m>0$. If  $u$ is a positive solution of  $(\ref{An1})$ in $\Gw$ and $z\in\prt\Gw$, we have that:   \smallskip

\nind 1- If there exists $\ge>0$ such that $|\nabla u|^q\in L^1_\gr(B_\ge(z)\cap \Gw)$, then $u^p\in L^1_\gr(B_\ge(z)\cap \Gw)$ and $u$ admits a boundary trace on $\prt\Gw\cap B_\ge(z)$ which is a nonnegative Radon measure.\smallskip

\nind 2- If there exists $\ge>0$ such that  $u^p\in L^1_\gr(B_\ge(z)\cap \Gw)$, then $u$ admits a boundary trace on $\prt\Gw\cap B_\ge(z)$ which is a nonnegative outer regular Borel measure, not necessarily bounded.\smallskip
\es

The last assertion shows how delicate is the construction of solutions with unbounded boundary trace.  We give a few examples with one point blow-up on the boundary. In particular we prove that when $0\in \prt\Gw$, $p>1$ and $\frac{2p}{p+1}<q<\frac{N+1}{N}$ there exist positive solutions $u$ of $(\ref{An1})$ in $\Gw$ (or $\Gw\setminus K$ where $K$ is compact), vanishing on $\prt\Gw\setminus\{0\}$ satisfying 
\begin{equation}\label{D03}
\BA{lll}\displaystyle
u(x)=|x|^{-\frac{2-q}{q-1}}\chi(\tfrac{x}{|x|})(1+o(1))\quad\text{as }\,x\to 0,
\EA\ee
for some positive function $\chi$ defined on the $S^{N-1}_+$. Such solutions have boundary trace $Tr_{\prt\Gw}(u)=(\{0\},0)$. \\
The existence of a boundary trace in the case $q> \frac{2p}{p+1}$ for {\it any} positive solution of $(\ref{An1})$ remains an open problem. \medskip

\nind{\it Acknowledgements}. The authors are grateful to the anonymous referees for the careful reading of the manuscript and their suggestions which lead to a substancial improvement of its presentation.

\mysection{Solutions with a boundary isolated singularity}
\subsection{A priori estimates}
In this section $\Gw$ denotes a bounded smooth domain of $\BBR^N$ such that $0\in\prt\Gw$. We prove an a priori estimate for positive solutions of  $(\ref{An1})$ vanishing on $\prt\Gw\setminus\{0\}$.
 \bprop{BSp} Suppose $1<p<\frac{N+2}{N-2}$ and either $1<q<\frac{2p}{p+1}$ and $m>0$, or $q=\frac{2p}{p+1}$ and $0<m<\ge_0$ for some $\ge_0>0$ depending on $N$ and $p$. Then there exists a constant $c=c(N,p,\Gw)>0$ such that if $u$ is a positive solution 
of $(\ref{An4})$, it satisfies
 \begin{equation}\label{Y2}
\BA{lll}
u(x)\leq c|x|^{-\ga}\qquad\text{for all }x\in\overline\Gw\setminus\{0\}.
\EA\ee
\es

The proof needs a series of intermediate results involving the Polacik et al. method \cite{PQS}, a result of Montoro \cite{Mon} and a previous Liouville 
theorem proved in \cite{BVGHV1}. We first recall the doubling lemma.
\blemma{PQS2019} Let $(X,d)$ be a complete metric space, $\Gg\varsubsetneq X$ and $\gg:X\setminus\Gg\mapsto (0,\infty)$. Assume that $\gg$ is bounded on all compact subsets of $X\setminus\Gg$. Given $k>0$, let $y\in X\setminus\Gg$ such that
$$\gg(y)\dist(y,\Gg)>2k.
$$
Then there exists $x\in X\setminus\Gg$ such that \\
1- $\gg(x)\dist(x,\Gg)>2k$,\\
2- $\gg(x)\geq\gg(y)$,\\
3- $2\gg(x)\geq\gg(z)$, for all $z\in B_{\frac{k}{\gg(x)}}(x)$.
\es
The next result is an extension of \cite[Proposition 5.1]{BVPV}. 
\blemma{BVPV2011} Suppose $1<p<\frac{N+2}{N-2}$ and either $1<q<\frac{2p}{p+1}$ and $m>0$, or $q=\frac{2p}{p+1}$ and $0<m<\ge_0$ for some $\ge_0>0$ depending on $N$ and $p$. Let $0<r<\frac{1}{2}\text{diam }\Gw$. There exists a constant $c>0$ depending on $p$, $m$, $q$ and $\Gw$ such that any  function $u$ verifying 
 \begin{equation}\label{Y3}
\BA{lll}
-\Gd u+m|\nabla u|^{q}=u^p\quad&\text{in }\Gw\cap (B_{2r}\setminus \overline B_r)\\
\phantom{-\Gd u+m|\nabla u|^{q}}
u\geq 0&\text{in }\Gw\cap (B_{2r}\setminus\overline  B_r)\\
\phantom{-\Gd u+m|\nabla u|^{q}}
u=0&\text{in }\prt\Gw\cap (B_{2r}\setminus\overline  B_r),
\EA\ee
satisfies
 \begin{equation}\label{Y4}
\BA{lll}
u(x)\leq c\left(\dist(x,\Gg_r)\right)^{-\ga}\quad\text{for all $x$ in } \Gw\cap (B_{2r}\setminus\overline  B_r),
\EA\ee
where $\Gg_r=\overline\Gw\cap\left(\prt B_{2r}\cup \prt B_r\right)$.
\es
\Proof We proceed by contradiction. For every $k\geq 1$ there exist $0<r_k<\frac{1}{2}\text{diam }\Gw$, a solution $u_k$ of $(\ref{Y3})$ with $r=r_k$ and 
$y_k\in \Gw\cap (B_{2r_k}\setminus\overline  B_{r_k})$ such that 
$$u_k(y_k)\geq (2k)^{\ga}\left(\dist(x,\Gg_{r_k})\right)^{-\ga}.
$$
It follows from \rlemma{PQS2019} applied with 
$$X=\overline\Gw\cap (\overline B_{2r_k}\setminus  B_{r_k})\quad\text{and }\; \gg=u_k^{\frac 1\ga},
$$
that there exists $x_k\in X\setminus\Gg_k$ such that 
 \begin{equation}\label{Y5}
\BA{lll}
\!\!\!\!\!(i)\; &u_k(x_k)\geq (2k)^{\ga}\left(\dist(x_k,\Gg_{r_k})\right)^{-\ga},\\
\!\!\!\!\!(ii)\; &u_k(x_k)\geq u_k(y_k),\\
\!\!\!\!\!(iii)\; &2^{\ga}u_k(x_k)\geq u_k(z),\;\text{for all }z\in B_{R_k}(x_k)\cap\Gw\;\text{with }R_k=k(u_k(x_k))^{-\frac{1}{\ga}}.
\EA\ee
Since (i) holds, $R_k<\frac 12\dist(x_k,\Gg_{r_k})$, hence 
$$B_{R_k}(x_k)\cap\Gg_{r_k}=\emptyset. $$
Since $\dist(x_k,\Gg_{r_k})\leq \frac 12r_k<\frac 14\text{diam }\Gw$, we also have from (i),
$$u_k(x_k)\geq\left(\myfrac{8k}{\text{diam }\Gw}\right)^{\ga}\to\infty\text{ as }k\to\infty.
$$
Next we set 
$$\BA {lll}
t_k=(u_k(x_k))^{-\frac 1\ga}\,,\,\; D_k=\left\{\xi\in\BBR^N:|\xi|\leq k\text{ and }x_k+t_k\xi\in\Gw\right\},
\EA$$
and 
$$v_k(\xi)=t_k^\ga u_k(x_k+t_k\xi)\quad\text{for all }\xi\in D_k.
$$
Then $v_k$ is positive in $D_k$ and satisfies 
 \begin{equation}\label{Y6}
\BA{lll}
-\Gd v_k+mt_k^{\frac{2p-(p+1)q}{p-1}}|\nabla v_k|^{q}=v_k^p\quad&\text{in } D_k\\[1mm]
\phantom{-}
0\leq v_k\leq 2^\ga\;&\text{ in  }D_k\\[1mm]
\phantom{-}v_k(0)=1.
\EA\ee
We encounter the following dichotomy:\smallskip

\nind 
(A) Either for every $a>0$ there exists ${k_a}\geq 1$ such that for $k\geq k_a$ $B_{at_{k}}(x_{k})\cap\prt\Gw=\emptyset$. The sequence $\{v_{k}\}$ is locally uniformly bounded  
in $\BBR^N$. Since $q\leq 2$, standard a priori estimates in elliptic equations imply that $\{v_{k}\}$ is eventually uniformly bounded in the $C^{2,\gt}$ local topology in $\BBR^N$ with $0<\gt<1$. Up to a subsequence still denoted by $\{v_{k}\}$ it converges locally in $C^2(\BBR^N)$ to a positive function $v$ which satisfies 
$v(0)=1$, $0\leq v\leq 2^\ga$ and either
 \begin{equation}\label{Y7}
\BA{lll}
-\Gd v=v^p\quad&\text{in } \BBR^N
\EA\ee
if $1<q<\frac{2p}{p+1}$, or 
 \begin{equation}\label{Y8}
\BA{lll}
-\Gd v+m|\nabla v|^{q}=v^p\quad&\text{in } \BBR^N
\EA\ee
if $q=\frac{2p}{p+1}$. In the first case it is proved in \cite{GS80} that such a solution cannot exist. If  $q=\frac{2p}{p+1}$ it is proved in \cite[Theorem E]{BVGHV1} that there exists $\ge_0>0$ depending on $N,p$ such that if $|m|\leq\ge_0$ no such solution exists. Therefore if situation (A) occurs we obtain a contradiction. \smallskip

\nind (B) Or there exists some $a_0>0$ such that $B_{a_0t_k}(x_k)\cap\prt\Gw\neq \emptyset$ for all $k\in\BBN^*$. Let $x'_k\in \prt\Gw$  minimizing the distance from $x_k$ and $\prt\Gw$, then 
$|x_k-x'_k|\leq a_0t_k$. Since the function $v_k$ is bounded in $D_k$ and vanishes on $\prt\Gw$ which is smooth,  it remains locally bounded in $W^{2,s}(\overline \Gw\cap B_{a'_0t_k}(x_k))$ for all $a'_0<a_0$ and all $s<\infty$, thus $\nabla v_k$ remains locally bounded therein. Then either $|x_k-x'_k|\geq a'_0t_k$ or $|x_k-x'_k|< a'_0t_k$. In this case we set $\xi'_k=t_k^{-1}(x'_k-x_k)$ and use the fact that  $v_k(\xi'_k)=0$ and 
$v_k(0)=1$ combined with the uniform bound on $\nabla v_k$ to infer that $|\xi'_k|\geq a_1$ for some $0<a_1<a_0$ independent of $k$ which implies that $|x_k-x'_k|\geq a_1t_k$. Up to a subsequence, we can assume that $t_k^{-1}x_k\to x_0$, $t_k^{-1}x'_k\to x'_0$ and that 
$D_k\to H$ where $H\sim \BBR_+^N$ is the half-space passing through $x'_0$ with normal inward unit vector ${\bf e}_{_N}$, and $x_0-x'_0=a{\bf e}_{_N}$ with $a_1\leq a\leq a_0$. Let $\tilde H\sim\BBR^N$ be the union of $H$ and its reflection through $\prt H$. Performing the reflection of $v_k$ through $\prt \left(t_k^{-1}\Gw\right)$ (see \cite[Lemma 3.3.2]{Vebook}) we deduce that the function $\tilde v_k$ which coincides with $v_k$ in $\prt \left(t_k^{-1}\Gw\right)$ and with its odd reflection in the image by reflection of the set $\prt \left(t_k^{-1}\Gw\right)$ 
vanishes on $\left(t_k^{-1}\Gw\right)$ and converges locally in $C^2(\BBR^N)$ to a positive function $\tilde v$ defined in $\tilde H$, bounded therein, vanishing on 
$\prt H$ and positive in $H$ and the function $v=\tilde v\lfloor_{H}$ is nonnegative and $v(x_0)=1$. If $q<\frac{2p}{p+1}$,
 $v$ satisfies 
 \begin{equation}\label{Y9}
\BA{lll}
-\Gd v=v^p\quad&\text{in } H\\
\phantom{-\Gd}v=0&\text{in } \prt H.
\EA\ee
By \cite{GS80} such a function cannot exist. If $q=\frac{2p}{p+1}$, the function $v$ satisfies 
 \begin{equation}\label{Y10}
\BA{lll}
-\Gd v+m|\nabla v|^{\frac{2p}{p+1}}=v^p\quad&\text{in } H\\
\phantom{-\Gd+m|\nabla v|^{\frac{2p}{p+1}}}v=0&\text{in } \prt H.
\EA\ee
Since it is positive, bounded and $\nabla v$ is also bounded, it follows from \cite{Mon} that $v$ is nondecreasing in the variable $x_{_N}$. 
But by  \cite[Theorem E]{BVGHV1}, the function $v$ satisfies 
 \begin{equation}\label{Y10-1}
\BA{lll}
v(x)\leq 2^\frac{2}{p-1}x_{_N}^{-\frac{2}{p-1}},
\EA\ee
which is impossible because $x_{_N}\mapsto v(., x_{_N})$ is nondecreasing. This ends the proof.\\$\phantom{--------}$ {\hspace{10mm}\hfill $\square$ \medskip

\nind{\it Proof of \rprop{BSp}}. We use \cite[Lemma 4.4]{BVPV} with $r=\frac{2|x|}{3}$ in $\Gw\cap\left(B_{\frac{4|x|}{3}}\setminus B_{\frac{2|x|}{3}}\right)$. $\phantom{----}$ {\hspace{10mm}\hfill $\square$} \medskip

\nind{\it Proof of \rth{BS1}}. It follows from standard regularity results and scaling techniques (see e.g. \cite[Lemma 3.3.2]{Vebook}). $\phantom{-oooooooo}$\qeda\medskip
\subsection{Removability}

\bth{BIS1} Let $\Gw$ be a bounded smooth domain such that $0\in\prt\Gw$, $\frac{N+1}{N-1}<p<\frac {N+2}{N-2}$ and either $1<q<\frac{2p}{p+1}$  and $m >0$, or $q=\frac{2p}{p+1}$ and $0<m < \ge^*$. If $u\in C^1(\overline\Gw\setminus\{0\})\cap C^2(\Gw)$ is a positive solution of 
$(\ref{An1})$ in $\Gw$ vanishing on $\prt\Gw\setminus \{0\}$, then $u\in L_\gr^p(\Gw)$, $\nabla u\in L_\gr^q(\Gw)$ and the equation $(\ref{An1})$ holds in the sense that 
 \begin{equation}\label{X4'}
\BA{lll}
\myint{\Gw}{}\left(-u\Gd\gz+\left(u^p-m|\nabla u|^q\right)\gz\right) dx=0\quad\text{for all }\,\gz\in \BBX(\Gw).
\EA\ee
\es
\Proof Under the assumptions on $p$ and $q$, there holds
  \begin{equation}\label{C1}
u(x)\leq c|x|^{-\ga}\qquad\text{for all }0<|x|\leq R. 
\ee
Since $p>\frac{N+1}{N-1}>\frac{N+2}{N}$, the function $u$ belongs to $L^1(\Gw)\cap L_\gr^p(\Gw)$. It follows by \rth{BS1} that, 
  \begin{equation}\label{C2}
|\nabla u(x)|\leq c|x|^{-\ga-1}\qquad\text{for all }0<|x|\leq \frac R2. 
\ee
Since $q\leq\frac{2p}{p+1}$ we have that  $|\nabla u(x)|^q\leq  c^q|x|^{-\ga-2}$. Hence $\nabla u$ belongs to $L^q_{\gr}(\Gw)$. Finally, let $\{\gz_n\}$ be a sequence of smooth functions such that $0\leq\gz_n\leq 1$, $\gz_n(x)=0$ if $|x|\leq n^{-1}$,  $\gz_n(x)=1$ if $|x|\geq 2n^{-1}$ with $|\nabla \gz_n(x)|\leq cn$ and $|\Gd\gz_n(x)|\leq cn^2$. Let $\gf\in \BBX(\Gw)$.  We have that
   \begin{equation}\label{C3}\BA{lll}
\myint{\Gw}{} u\Gd(\gz_n\gf) dx=\myint{\Gw}{}\gz_n u\Gd\gf dx+\myint{\Gw}{}\gf u\Gd\gz_n  dx+2\myint{\Gw}{}u\nabla \gf.\nabla\gz_n dx\\[2mm]
\phantom{\myint{\Gw}{} u\Gd(\gz_n\gf) dx}=I(n)+II(n)+III(n).
\EA\ee
Clearly 
$$I(n)\to \myint{\Gw}{} u\Gd\gf dx\quad\text {as }n\to\infty.
$$
If $\gf\in \BBX(\Gw)$, $\gr^{-1}\gf$ is bounded in $\Gw$, hence
$$\BA{lll}
|II(n)|\leq c_1 n^2\norm{\gr^{-1}\gf}_{L^\infty}\left(\myint{n^{-1}\leq |x|\leq 2n^{-1}}{}u^p\gr dx\right)^\frac{1}{p}\left(\myint{n^{-1}\leq |x|\leq 2n^{-1}}{} \gr dx\right)^\frac{1}{p'}\\[2mm]
\phantom{|II(n)|} \leq c'_1 n^{2-\frac{N+1}{p'}}\norm{\gr^{-1}\gf}_{L^\infty}\left(\myint{n^{-1}\leq |x|\leq 2n^{-1}}{}u^p\gr dx\right)^\frac{1}{p}.
\EA$$
Since $p>\frac{N+1}{N-1}$, $p'<\frac{N+1}{2}$, hence $|II(n)|\to 0$ when $n\to\infty$. For the last term, we have from 
\rth{BS1},
$$\BA{lll}|III(n)|\leq c_2n\norm{\nabla\gf}_{L^\infty}\myint{n^{-1}\leq |x|\leq 2n^{-1}}{}u dx
\\\phantom{|III(n)|} 
\leq c'_2n^{}\norm{\nabla\gf}_{L^\infty}\myint{n^{-1}\leq |x|\leq 2n^{-1}}{}|x|^{-\ga-1}\gr dx\\[2mm]\phantom{|III(n)|} 

\leq c''_2n^{\ga+1-N}\norm{\nabla\gf}_{L^\infty}.
\EA$$
Since $p>\frac{N+1}{N-1}$, $\ga+1-N<0$, we deduce that $|III(n)|\to 0$ when $n\to\infty$.
Therefore there holds 
$$\myint{\Gw}{}\left(u^p-m |\nabla u|^q\right) \gf\gz_ndx\to \myint{\Gw}{}\left(u^p-m |\nabla u|^q\right) \gf dx,
$$
we obtain the claim.\qeda\medskip

\subsection{Proof of \rcor{BSS1}}
The proof is an easy but technical adaptation of the computations in \cite[Theorems 1.1, 1.2]{BVPV} and \cite[Theorem 3.25]{NGT-LV}, but for the sake of completeness, we briefly recall its technique. Since $\Gw$ is in normal position with respect to $0$ there exist a bounded open neighborhood $G$ of $0$ and a smooth function $\gf: G\cap\prt\BBR_+^N\mapsto\BBR$ such that 
$$G\cap\prt\Gw=\left\{x=(x',x_{_N}):x'\in G\cap\prt\BBR_+^N\text{ and }x_{_N}=\phi(x')\right\}.
$$
Furthermore $\phi(x')=0(|x'|^2)$, $\nabla \phi(x')=0(|x'|)$ and $|D^2\phi(x')|\leq c$ if $x'\in G\cap\prt\BBR_+^N$. If $u$ satisfies $(\ref{An1})$, we 
denote
$$u(x)=\tilde u(y)\text{ with }y_i=x_i\text{ when }1\leq i\leq N-1\text{ and }y_{_N}=x_{_N}-\phi(x').$$
If we set $r=|y|$, $s=y/r$, $t=\ln r$ and $v(t,s)=r^\ga\tilde u(r,s)$, then $v$ is bounded in $C^{2}((-\infty,T_0]\ti S^{N-1}_+)$ and vanishes on 
$(-\infty,T_0]\ti \prt S^{N-1}_+)$. Using the computations in \cite[Theorem 3.25]{NGT-LV} and \cite[Lemma 6.1]{BVPV}, it satisfies, with 
${\bf n}=\frac{y}{|y|}$,
\bel{X1}\BA {lll}
(1+\ge_1)v_{tt}+\Gd' v-(N-2+2\ga+\ge_2)v_t+\left(\ga(N-2-\ga)+\ge_3\right)v\\[2mm]\phantom{-------}
+\Gd'v+\nabla' v.\overrightarrow\ge_4+\nabla' v_t.\overrightarrow\ge_5+\nabla'(\nabla' v.{\bf e}_{_N}).\overrightarrow\ge_6+v^p
\\[2mm]\phantom{-}
-me^{\frac{2p-q(p+1)}{p-1}t}\left[\phantom{b^{p^p}}\!\!\!\!\!\!\!\!\left(v_t-\ga v\right){\bf n}+\nabla' v+\left(\left(v_t-\ga v\right){\bf n}+\nabla' v.{\bf e}_{_N}\right).\overrightarrow\ge_7\right]^{q}
=0,
\EA\ee
where $\CB:=\{{\bf e}_{_1},...,{\bf e}_{_N}\}$ denotes the canonical orthogonal basis in $\BBR^N$. The functions $\ge_j$ (or $\overrightarrow\ge_j$) are uniformly continuous and bounded for $j=1,...,7$ and there holds
\bel{X2}\BA {lll}
|\ge_j(t,.)|\leq ce^t\quad&\text{for }j=1,...,7\\
|\ge_{jt}(t,.)|+|\nabla '\ge_j(t,.)|\leq ce^t\quad&\text{for }j=1,5,6,7.
\EA\ee
By since $v$,$v_t$ and $\nabla'v$ are uniformly bounded, we infer by standard regularity results (see e.g. \cite{GT}) the following uniform estimate,
\bel{X3}\BA {lll}
\norm{v(t,.)}_{C^{2,\gt}(\overline{S^{N-1}_+})}+\norm{v(t,.)_t}_{C^{1,\gt}(\overline{S^{N-1}_+})}+\norm{v_{tt}(t,.)}_{C^{0,\gt}(\overline{S^{N-1}_+})}
\leq c
\EA\ee
for any $t\leq T_0$, for some $c>0$ and $\gt\in (0,1)$. Hence the limit set at $-\infty$ of the trajectory $\{v(t,.)\}_{t\leq T_0}$ in $C^{2}(\overline{S^{N-1}_+})$ is a connected non-empty compact subset of $\left\{\gw\in C^2(\overline{S^{N-1}_+}):\gw\lfloor_{\prt S^{N-1}_+}=0\right\}$. Next we write $(\ref{X1})$ under the form
\bel{X4}\BA {lll}
v_{tt}+\Gd' v-(N-2+2\ga)v_t+\ga(N-2-\ga)v+v^p=e^{\gth t}\Gth,\EA\ee
where $\Gth$ is bounded and $\gth=\min\left\{1,\frac{2p-q(p+1)}{p-1}\right\}$. Since $N-2+2\ga\neq 0$, the standard energy method (multiplication by 
$v_t$) yields 
$$\myint{-\infty}{T_0}\myint{S^{N-1}_+}{}(v_t^2+v^2_{tt})dSdt<\infty.
$$
Since $v_{t}$ and $v_{tt}$ are uniformly continuous, the above integrability condition yields
\bel{X5}\BA {lll}
\displaystyle \lim_{t\to-\infty}\left(\norm {v_t(t,.)}_{L^2(S^{N-1}_+)}+\norm {v_{tt}(t,.)}_{L^2(S^{N-1}_+)}\right)=0.
\EA\ee
Therefore the limit set of the trajectory at $-\infty$ is a compact connected subset of nonnegative solutions of $(\ref{An10})$. This implies that 
either $v(t,.)$ converges to the unique positive solution $\psi$ of $(\ref{An10})$ in $C^{2}(\overline{S^{N-1}_+})$ or it converges to $0$. 
Note that the set of nonnegative solutions of $(\ref{An10})$ is reduced to $0$ when $1<p\leq\frac{N+1}{N-1}$. \smallskip

\nind If $\frac{N+1}{N-1}<p<\frac{N+2}{N-2}$ and $v(t,.)$ does not converge to $0$, then we have proved $(\ref{An9})$. If $v(t,.)$ converges to $0$, then the proof of \cite[Theorem 7.1]{BVPV} applies, the only difference being in the value of the term $H$ therein \cite[(7.3)]{BVPV} 
which is replaced by $e^{\gth t}\Gth$ defined above. The remaining of the argument can be easily adapted. \smallskip

\nind If $p=\frac{N+1}{N-1}$ then $v(t,.)$ converges to $0$. The adaptation of \cite[Theorem 9.1]{BVPV} is easy. We obtain that $u$ satisfies 
\bel{X6}\BA {lll}
u(x)\leq c|x|^{1-N}\left(\ln \frac{1}{|x|}\right)^{-\frac{N-1}{2}}\qquad\text{for all }x\in\Gw.
\EA\ee
The completion of the proof follows by the same perturbation method as in \cite[Lemma 9.1]{BVPV}, by decomposing the function $v(t,.)$ into $v(t,.)=v_1+v_2(t,.)$ where 
$v_1\in\ker(\Gd'+(N-1)I)$ and $v_2\in(\ker(\Gd'+(N-1)I))^{\perp}$. This yields
\bel{X7}\BA {lll}
\norm {v_1(t,.)}_{L^2(S^{N-1}_+)}\!\leq c(-t)^{-\frac{N-1}{2}}\text{ and }\,\norm {v_2(t,.)}_{L^2(S^{N-1}_+)}\!\leq ce^{\frac{\gth}{2}t}\;\text{ for }t\leq T_0.
\EA\ee
The function $w(t,s)=(-t)^{\frac{N-1}{2}}v(t,s)$ satisfies 
\bel{X8}\BA {lll}
w_{tt}-\left(N+\frac{N-1}{t}\right)w_t+\left(N-1+\frac{N^2-1}{4t^2}\right)w+\Gd'w\\[2mm]
\phantom{---------}-\myfrac{1}{t}\left(w^{\frac{N+1}{N-1}}-\frac{N(N-1)}{2}w\right)
=(-t)^{\frac{N-1}{2}}\Gth,
\EA\ee
where $\Gth$ is bounded. The proof given in \cite[Theorem 9.1]{BVPV} applies with almost no change, but for  some straightforward ones. The main step is to introduce 
$$z(t)=\myint{S^{N-1}_+}{}w(t,s)\phi_1(s)dS,
$$
and to prove that $z(t)$ admits a nonnegative  limit $\gl\geq 0$ when $t\to-\infty$. If  this limit is positive its value $\gl$ is given in the proof of  \cite[Theorem 1.3]{BVPV}. If this limit is zero, then 
$$\displaystyle\lim_{y\to 0}|y|^{N-1}\left(\ln\frac1{|y|}\right)^{\frac{N-1}{2}}\tilde u(y)=0,
$$
and the conclusion follows easily from the proof \cite[Theorem 7.2]{BVPV} (only the exponent in the perturbation term $H$ therein is changed).\smallskip

\nind If $1<p<\frac{N+1}{N-1}$, then $v(t,.)$ converges to $0$ and $(\ref{X4})$ can be written under the form 
\bel{X9}\BA {lll}
v_{tt}+\Gd' v-(N-2+2\ga)v_t+\ga(N-2-\ga+\ge(t))v+\Gd'v=0,
\EA\ee
where $\ge(t)\to 0$ when $t\to-\infty$. It is therefore a very standard but technical method of linearization \cite[Theorem 5.1]{GmVe} to obtain, first an exponential decay of $w(t,.)$ at $-\infty$, and then the convergence of $t\mapsto e^{(N-1-\ga)t}v(t,.)$ to $k\phi_1$ for some 
$k\geq 0$, and then to deduce the regularity of $u$ if $k=0$.\qeda
\mysection{Separable solutions}
\subsection{Separable solutions in $\BBR^N$}
\nind{\it Proof of \rprop{rad}}. Constant positive solutions of  $(\ref{A6})$ are any positive roots of
\bel{fi1}
\Gf(X):=X^{p-1}-m\ga^{\frac{2p}{p+1}}X^{\frac{p-1}{p+1}}-\ga(N-2-\ga)=0.
\ee
Set
\bel{fi2}\Gf(X)=\tilde \Gf(X^{\frac{p-1}{p+1}}),\ee
where 
\bel{fi3}\tilde \Gf(Y)=Y^{p+1}-m\ga^{\frac{2p}{p+1}}Y-\ga(N-2-\ga).\ee
Then 
$\tilde \Gf'(Y)=(p+1)Y^p-m\ga^{\frac{2p}{p+1}}$, hence if $m\leq 0$, $\tilde \Gf$ is increasing and if $m>0$, $\tilde \Gf$ is decreasing on $[0,Y_0)$ and increasing on $(Y_0,\infty)$ with 
\begin{equation}\label{A3}
Y_0=\left(\myfrac{m}{p+1}\right)^{\frac{1}{p}}\ga^{\frac{2}{p+1}}.
\ee
From now we always assume $m>0$. Then 
$$\BA {lll}\tilde\Gf(Y_0)=\left[N-p(N-2)-2p\left(\myfrac{m}{p+1}\right)^{\frac{p+1}{p}}\right]\myfrac{2}{(p-1)^2},
\EA$$
and 
$$\tilde\Gf(0)=-\ga(N-2-\ga)=\myfrac{2(N-2)}{(p-1)^2}\left(\myfrac{N}{N-2}-p\right).$$
Therefore, $\tilde\Gf(0)\leq 0$ if and only if $p\geq \frac{N}{N-2}$. In that case there exists a unique $X_{_m}>0$ such that  $\Gf(X_{_m})=0$. \\
When 
\begin{equation}\label{A'4}
0<\myfrac{N}{N-2}-p<\myfrac{2p}{N-2}\left(\myfrac{m}{p+1}\right)^{\frac{p+1}{p}},
\ee
then $\tilde\Gf(0)>0$ and $\tilde\Gf(Y_0)<0$, thus $\tilde \Gf$ admits two positive roots. The same property is shared by $\Phi$, hence there exist $X_{{j,m}}$, for $j=1,2$ such that $\Phi(X_{{j,m}})=0$ and $0<X_{{1,m}}<Y_0^{\frac{p+1}{p-1}}<X_{{2,m}}$. \\
When
\begin{equation}\label{A5}
0=\myfrac{2p}{N-2}\left(\myfrac{m}{p+1}\right)^{\frac{p+1}{p}}\Longleftrightarrow \left(\myfrac{m}{p+1}\right)^{\frac{p+1}{p}}=\myfrac{N-p(N-2)}{2p},
\ee
then $\tilde\Gf$ admits a unique positive root. Hence $\Gf>0$ on $\BBR_+\setminus\{X_{m^*}\}$ and vanishes at $X_{m^*}$, where 
 \begin{equation}\label{Ax6}
X_{m^*}=\left(\myfrac{m^*}{p+1}\right)^{\frac{p+1}{p'p-1)}}\ga^{\frac{2}{p-1}}\quad\text{with }\,m^*=(p+1)\left(\myfrac{N-p(N-2)}{2p}\right)^{\frac{p}{p+1}}.
\ee
If $0<m<m^*$, $\tilde \Gf$ and thus $\Gf$ are positive on $\BBR_+$, hence there exists no root to $\Gf$.
The proof of \rprop{rad} is complete.\qeda
\subsection{Separable solutions in $\BBR_+^N$}
If $u$ is a nonnegative separable solution of $(\ref{An1})$ in $\BBR^N_+$ which vanishes on $\prt\BBR^N_+\setminus\{0\}$, the function $\gw$ is a nonnegative solution of $(\ref{A9})$.\smallskip

\nind {\it Proof of \rth{Exist}}. If $1<\frac{2p}{p+1}<\frac{N+1}{N}$, equivalently $1<p<\frac{N+1}{N-1}$, it is proved in \cite[Theorem 3.21]{NGT-LV} that there exists a unique positive function $\eta:=\eta_m\in C^2(\overline{S_+^{N-1}})$ satisfying 
 \begin{equation}\label{Z4}
\BA{lll}
-\Gd' \eta+\ga(N-2-\ga)\eta+m\left(\ga^2\eta^2+|\nabla' \eta|^2\right)^{\frac{p}{p+1}}=0\quad&\text{in }\;S_+^{N-1}\\
\phantom{-\Gd' +\ga(N-2-\ga)\eta+m\left(\ga^2\eta^2+|\nabla' \eta|^2\right)^{\frac{p}{p+1}}}
\eta=0&\text{on }\;\prt S_+^{N-1}.
\EA\ee
By uniqueness, $\eta_m=m^{-\frac{p+1}{p-1}}\eta_1$, and by the maximum principle
 \begin{equation}\label{Z5}
\BA{lll}
m^{-\frac{p+1}{p-1}} \displaystyle\sup_{S^{N-1}_+}\eta_1= \displaystyle\sup_{S^{N-1}_+}\eta_m\leq \myfrac{1}{\ga}\left(\myfrac{\ga+2-N}{m}\right)^{\frac{p+1}{p-1}}.
\EA\ee
If $\overline \eta_m=\displaystyle\sup_{S^{N-1}_+}\eta_m$, then 
$$-m\ga^{\frac{2p}{p+1}}\overline \eta_m^{\frac{2p}{p+1}}-\ga(N-2-\ga)\overline \eta_m\geq 0.
$$
Hence $\Phi(\overline \eta_m)>0$, where $\Phi$ has been defined in $(\ref{fi1})$. Therefore  
 \begin{equation}\label{Z6}
\BA{llll}
(i)  &\text{either }\qquad &\overline \eta_m>X_{2,m}\quad&\text{(resp. $\overline \eta_{m^*}>  X_{m^*}$)},\quad\qquad\qquad\quad\quad\quad\\[2mm]
(ii) &\text{or }\qquad &\overline\eta_m<X_{1,m}\quad&\text{(resp. $\overline\eta_{m^*}<  X_{m^*}$)}.\\[2mm]
\EA\ee
For  $\ge\in (0,1)$, $\ge\eta_m$ is a subsolution of $(\ref{Z4})$, hence it is a subsolution of $(\ref{A9})$ too. For $\ge>0$ small enough it is smaller than $X_{2,m}$ (resp. $X_{m^*}$) and it belongs to $W^{1,\infty}_0(S^{N-1}_+)$.
By the result of Boccardo, Murat and Puel \cite{BMP} there exists a solution $\gw\in W^{1,2}_0(S^{N-1}_+)$ of  $(\ref{A9})$, and it satisfies 
 \begin{equation}\label{Z7}
\BA{lll}
\ge\eta_m<\gw\leq  X_{2,m}\quad\text{(resp. $\ge\eta_m<\gw\leq  X_{m^*}$)}.
\EA\ee

For proving the second assertion, 
we set $\gw=\phi^b$ for some $b>1$, then 
 \begin{equation}\label{Z8}
\BA{lll}
-\Gd'\phi-(b-1)\myfrac{|\nabla\gf|^2}{\gf}-\myfrac{\ga(\ga+2-N)}{b}\gf-\myfrac{1}{b}\phi^{1+b(p-1)}\\[4mm]
\phantom{---------}+\myfrac{m}{b}\phi^{\frac{(p-1)(b-1)}{p+1}}\left(\ga^2\gf^2+b^2|\nabla\gf|^2\right)^{\frac{p}{p+1}}=0.
\EA\ee
Since
 \begin{equation}\label{Z9}
\BA{lll}
\left(\ga^2\gf^2+b^2|\nabla\gf|^2\right)^{\frac{p}{p+1}}\leq \ga^{\frac{2p}{p+1}}\gf^{\frac{2p}{p+1}}+b^{\frac{2p}{p+1}}|\nabla\gf|^{\frac{2p}{p+1}},
\EA\ee
$(\ref{Z8})$ implies
 \begin{equation}\label{Z11}
\BA{lll}
-\Gd'\phi+\myfrac{m\ga^{\frac{2p}{p+1}}}{b}\gf^{1+b\frac{p-1}{p+1}}+mb^{\frac{p-1}{p+1}}\gf^{\frac{(b-1)(p-1)}{p+1}}|\nabla\gf|^{\frac{2p}{p+1}}\\[4mm
]\phantom{-----} \geq(b-1)\myfrac{|\nabla\gf|^2}{\gf}+\myfrac{1}{b}\phi^{1+b(p-1)}+\myfrac{\ga(\ga+2-N)}{b}\gf.
\EA\ee
For any $\gth>0$ we have by H\"older's inequality,
$$\BA {lll}mb^{\frac{p-1}{p+1}}\gf^{\frac{(b-1)(p-1)}{p+1}}|\nabla\gf|^{\frac{2p}{p+1}}\leq 
\myfrac{mpb^{\frac{p-1}{p+1}}}{(p+1)\gth^{\frac{p+1}{p}}}\myfrac{|\nabla\gf|^2}{\phi}+\myfrac{mb^{\frac{p-1}{p+1}}\gth^{p+1}}{p+1}\phi^{1+b(p-1)},
\EA$$
we deduce the inequality 
 \begin{equation}\label{Z12}
\BA{lll}
-\Gd'\phi\geq \left(b-1-\myfrac{mpb^{\frac{p-1}{p+1}}}{(p+1)\gth^{\frac{p+1}{p}}}\right)\myfrac{|\nabla\gf|^2}{\phi}+\myfrac{1}{b}\left(1-\myfrac{mb^{\frac{2p}{p+1}}\gth^{p+1}}{p+1}\right)\phi^{1+b(p-1)}\\[4mm]\phantom{------------------------}
+\myfrac{\ga(\ga+2-N)}{b}\gf.
\EA\ee
If the following two conditions are satisfied
 \begin{equation}\label{Z13}
\BA{lll}
(i)\qquad &b-1-\myfrac{mpb^{\frac{p-1}{p+1}}}{(p+1)\gth^{\frac{p+1}{p}}}\geq 0,\qquad  \qquad \qquad \qquad \qquad \qquad \qquad \qquad \\[5mm]
(ii)\qquad &1-\myfrac{mb^{\frac{2p}{p+1}}\gth^{p+1}}{p+1}\geq 0,
\EA\ee
we infer that there holds
 \begin{equation}\label{Z14}
\BA{lll}
(N-1)\myint{S^{N-1}_+}{}\phi\phi_1 dS> \myfrac{\ga(\ga+2-N)}{b}\myint{S^{N-1}_+}{}\phi\phi_1 dS,
\EA\ee
where $\phi_1$ denotes the first normalized and positive eigenfunction of $-\Gd'$ in $W^{1,2}_0(S^{N-1}_+)$, with corresponding eigenfunction $\gl_1=N-1$. Hence,  if $(\ref{Z13})$ is verified and there holds
 \begin{equation}\label{Z15}
\BA{lll}
N-1\leq \myfrac{\ga(\ga+2-N)}{b},
\EA\ee
there exists no positive solution. We proceed as follows for solving $(\ref{Z13})$-$(\ref{Z15})$. If $1<p<\frac{N+1}{N-1}$, then $\ga(\ga+2-N)>N-1$. We define 
$b_p>1$ by 
 \begin{equation}\label{Z16}
\BA{lll}
b_p=\myfrac{\ga(\ga+2-N)}{N-1}.
\EA\ee
For such $b=b_p$, the optimality is achieved in $(\ref{Z13})$ when $b_p-1=\frac{mpb_p^{\frac{p-1}{p+1}}}{(p+1)\gth^{\frac{p+1}{p}}}$ and 
$1=\frac{mb_p^{\frac{2p}{p+1}}\gth^{p+1}}{p+1}$. This gives an implicit maximal value of $ m_p$ through the relation
 \begin{equation}\label{Z17}m_p=\myfrac{(p+1)(b_p-1)\gth^{\frac{p+1}{p}}}{pb_p^{\frac{p-1}{p+1}}}=\myfrac{p+1}{b_p^{\frac{2p}{p+1}}\gth^{p+1}}.
\ee
Then the value of the corresponding $\gth:=\gth_p$ is expressed by 
$$\gth_p=\myfrac{p}{b_p(b_p-1)},
$$
and we infer
 \begin{equation}\label{Z18}m_p=\myfrac{p+1}{b_p^{\frac{2p}{p+1}}\gth_p^{p+1}}=\myfrac{(p+1)}{p^{p+1}}(b_p-1)^{p+1}b_p^{\frac{p^2+1}{p+1}}.
\ee
Hence if $m\leq m_p$ problem $(\ref{Z9})$ admits no positive solution.\qeda\medskip

\nind \Remark The case $p\geq \frac{N+1}{N-1}$ is open. It can be noticed that the constant solution $X_m$ obtained in \rprop{rad}-(i) cannot be used as a supersolution for solving problem $(\ref{A9})$ as it is done in \rth{Exist}. If $\gw$ is a positive solution of $(\ref{A9})$ and $\overline\gw$ is it maximal value, then
$$-\Gd\overline\gw=\overline\gw\Phi(\overline\gw).
$$
Hence $\Phi(\overline\gw)\geq 0$ which implies that $\overline\gw>X_m$. 
\mysection{Boundary data measures}
\subsection{Sufficient conditions}

We associate to  $(\ref{An2})$ the following two problems
\begin{equation}\label{D2}
\BA{lll}
-\Gd v+m|\nabla v|^{q}=0&\qquad\text{in }\Gw\\
\phantom{-\Gd +m|\nabla v|^{q}}
v=\gm&\qquad\text{in }\prt\Gw,
\EA\ee
and 
\begin{equation}\label{D3}
\BA{lll}
-\Gd w-w^p=0&\qquad\text{in }\Gw\\
\phantom{-\Gd-w^p}
w=\gm&\qquad\text{in }\prt\Gw.
\EA\ee
Problem $(\ref{D2})$ has been solved in the case $1<q<\frac{N+1}{N}$ in \cite{NGT-LV}. There, it is proved that for any nonnegative bounded measure $\gm$ on $\prt\Gw$ there exists a weak solution $v_\gm$ to   $(\ref{D2})$. Furthermore the correspondance $\gm\mapsto v_\gm$ is sequentially stable. When $ \frac{N+1}{N}\leq q<2$ it is proved in  \cite[Theorem 1.6]{BVHNV} that if a measure $\gm$ satisfies 
\begin{equation}\label{D4}
\BA{lll}
|\gm|(K)\leq C_1Cap^{\prt\Gw}_{\frac{2-q}{q},q'}(K)\quad\text{for any compact set } K\subset\prt\Gw,
\EA\ee
then there exists $\ge_0>0$ such that for any $0<\ge\leq\ge_0$ there exists a solution $v_{\ge\gm}$ to $(\ref{D2})$ (i.e. with $\gm$ replaced by $\ge\gm$).\\
Problem $(\ref{D3})$ has been solved in the case $1<p<\frac{N+1}{N-1}$ in \cite{BV-Vi} where it is proved that for any nonnegative measure $\gm$ there exists $\ge_1>0$ such that for any $0<\ge\leq\ge_1$ there exists a positive solution $w:=w_{\ge\gm}$ to $(\ref{D3})$ provided $\gm$ is replaced by $\ge\gm$. In the supercritical case $p\geq \frac{N+1}{N-1}$ it is shown in \cite[Theorem 1.6]{BVHNV} that if a positive measure $\gm$ satisfies 
\begin{equation}\label{D5}
\BA{lll}
\gm(K)\leq C_2Cap^{\prt\Gw}_{\frac{2}{p},p'}(K)\quad\text{for any compact set } K\subset\prt\Gw,
\EA\ee
then existence of a positive solution $w_{\ge\gm}$ to problem $(\ref{D3})$ holds with $\gm$ replaced by $\ge\gm$, under the condition $0<\ge\leq \ge_2$, for some $\ge_2>0$ depending on $\gm$. \medskip

\nind{\it Proof of \rth{meas}}. We assume that $(\ref{D6})$ holds and we set $\ge_3=\min\{\ge_0,\ge_1,\ge_2\}$, take $\ge\leq\ge_3$ and for the sake of clarity, replace 
$\ge\gm$ by $\gm$. We denote by $v_{\gm}$ and $w_{\gm}$ the solutions of $(\ref{D2})$ and $(\ref{D3})$ respectively with boundary data $\gm$. Since there holds
$$v_{\gm}\leq\BBP_\Gw[\gm]\leq w_{\gm},
$$
and $v_{\gm}$ is a subsolution of $(\ref{An1})$ and $w_{\gm}$ a supersolution in $\Gw$, it follows from \cite[Theorem 1.4.6]{Vebook} that there exists a solution $u$ to $(\ref{An1})$ such that $v_{\gm}\leq u\leq w_{\gm}$. This implies that $u\in L^1(\Gw)$ and $u^p\in L^1_\gr(\Gw)$. Because $v$ and $w$ satisfy 
\begin{equation}\label{D7}
\BA{lll}
\displaystyle
\lim_{\gd\to 0}\myint{\{\gr(x)=\gd\}}{}v ZdS=\lim_{\gd\to 0}\myint{\{\gr(x)=\gd\}}{}w ZdS=\myint{\prt\Gw}{}Z d\gm
\EA\ee
for any $Z\in C(\overline\Gw)$, it follows that 
\begin{equation}\label{D8}
\BA{lll}\displaystyle
\lim_{\gd\to 0}\myint{\{\gr(x)=\gd\}}{}u ZdS=\myint{\prt\Gw}{}Z d\gm.
\EA\ee
Let $\phi_\gd$ be the first eigenfunction of $-\Gd$ in $W^{1,2}_0(\Gw'_\gd)$ ($\Gw'_\gd$ is defined in $(\ref{D22})$), normalized by $0\leq\phi_\gd\leq 1=\max\{\phi_\gd(x) :x\in \Gw'_\gd\}$ and $\gl_\gd$ the eigenvalue. Then 
\begin{equation}\label{D9-3}m\myint{\Gw'_\gd}{}|\nabla u|^{q}\phi_\gd dx=
\myint{\Gw'_\gd}{}\left(u^p-\gl_\gd u\right)\phi_\gd dx-\myint{\Gs_\gd}{}\myfrac{\prt\phi_\gd}{\prt{\bf n}}u(x) dS.
\ee
Because $\phi_\gd\to \phi :=\phi_0$ and $\gl_\gd\to \gl :=\gl_0$, and the left-hand side of $(\ref{D9-3})$ is convergent, it follows by Fatou's lemma that 
$$
m\myint{\Gw}{}|\nabla u|^{q}\phi dx\leq
\myint{\Gw}{}\left(u^p-\gl u\right)\phi dx-\myint{\prt\Gw}{}\myfrac{\prt\phi}{\prt{\bf n}}d\gm.
$$
Hence $\nabla u\in L_\gr^q(\Gw)$, thus $(\ref{D1-1})$ holds and this ends the proof. \qeda\medskip

In several cases the sufficient condition can be weakened either by comparison between capacities or because one at least of the two exponents $p$ or $q$ is subcritical. \medskip

\nind {\it Proof of \rcor{cormeas1}}. It follows easily from \cite {BV-Vi}, \cite {NGT-LV} and the previous theorem.
$\phantom{------------}$\qeda
\medskip

\nind {\it Proof of \rcor{cormeas2}}  1- As in the proof of \cite[Corollary 1.5]{BVGHV3}, we have from \cite[Theorem 5.5.1]{AdHe} 
$$Cap^{\prt\Gw}_{\frac{2}{p},p'}(K)\leq c^*Cap^{\prt\Gw}_{\frac{2-q}{q},q'}(K).
$$
It implies the following inequality
$$\BA{lll}\gm(K)\leq C_3Cap^{\prt\Gw}_{\frac{2}{p},p'}(K)=C_3\min\left\{Cap^{\prt\Gw}_{\frac{2}{p},p'}(K),c^*Cap^{\prt\Gw}_{\frac{2-q}{q},q'}(K)\right\}\\[2mm]
\phantom{\gm(K)}
\leq C_3(1+c^*)min\left\{Cap^{\prt\Gw}_{\frac{2}{p},p'}(K),Cap^{\prt\Gw}_{\frac{2-q}{q},q'}(K)\right\}.
\EA$$
2- Similarly, as in the proof of \cite[Corollary 1.4]{BVGHV3}, we have from \cite[Theorem 5.5.1]{AdHe} 
$$Cap^{\prt\Gw}_{\frac{2-q}{q},q'}(K)\leq c^{**}Cap^{\prt\Gw}_{\frac{2}{p},p'}(K),
$$
therefore 
$$\BA{lll}\gm(K)\leq C_4Cap^{\prt\Gw}_{\frac{2-q}{q},q'}(K)=C_4\min\left\{c^{**}Cap^{\prt\Gw}_{\frac{2}{p},p'}(K),Cap^{\prt\Gw}_{\frac{2-q}{q},q'}(K)\right\}\\[2mm]
\phantom{\gm(K)}
\leq C_4(1+c^{**})min\left\{Cap^{\prt\Gw}_{\frac{2}{p},p'}(K),Cap^{\prt\Gw}_{\frac{2-q}{q},q'}(K)\right\}.
\EA$$
This completes the proof. \qeda

\subsection{Necessary conditions}
{\it Proof of \rth{meas-nec}}. \\
{\it Step 1: proof of $(\ref{D6-1})$}. We follow the notations of the proof of \cite[Theorem 4-5]{NGT-LV}. Let $\eta\in C^2(\prt\Gw)$ be a nonnegative function  with value $1$ in a neighborhood $\CU$ of the compact set $K$, and $\gz=(\BBP_\Gw[\eta])^{2q'}\phi$. Then we have
$$\myint{\Gw}{}\left(|\nabla u|^q\gz-u\Gd\gz\right)dx=\myint{\Gw}{}u^p\gz dx-\myint{\prt\Gw}{}\myfrac{\prt\gz}{\prt{\bf n}}dS\geq -\myint{\prt\Gw}{}\myfrac{\prt\gz}{\prt{\bf n}}dS.
$$
Since $\eta=1$ on $K$, there holds by Hopf boundary lemma
$$-\myint{\prt\Gw}{}\myfrac{\prt\gz}{\prt{\bf n}}dS\geq c_1\gm(K).
$$
The same computation as in \cite[Theorem 4-5]{NGT-LV} yields, with $\gl=\gl_1(\Gw)$, 
\begin{equation}\label{D6-3}c_1\gm(K)\leq 
\myint{\Gw}{}\left(|\nabla u|^q+\gl u\right)\gz dx+c_2\left(1+\norm{\nabla u}_{L^q_\gr}^q\right)^{\frac 1q}
\norm\eta_{W^{\frac{2-q}{q},q'}}.
\ee
Since $Cap^{\prt\Gw}_{\frac{2-q}{q},q'}(K)=0$, there exists a sequence $\{\eta_n\}\subset C^2(\prt\Gw)$ satisfying 
$0\leq \eta_n\leq 1$ and $\eta_n=1$ in a neighborhood of $K$, such that $\norm{\eta_n}_{W^{\frac{2-q}{q},q'}}\to 0$ as 
$n\to\infty$; thus $\eta_n\to 0$ in $L^1(\prt\Gw)$ and $\gz_n:=(\BBP_\Gw[\eta_n])^{2q'}\phi\to 0$ a.e. in $\Gw$. 
This implies that the right-hand side of $(\ref{D6-3})$ with $\eta$ replaced by $\eta_n$ tends to $0$ as $n\to\infty$ and thus $\gm(K)=0$.\smallskip

\nind {\it Step 2: proof of $(\ref{D6-2})$}. We recall that a positive lifting is a mapping $\eta\mapsto R[\eta]$ from $C^2(\prt\Gw)$ to  $C^2(\overline \Gw)$ satisfying 
$$R[\eta]\lfloor_{\prt\Gw}=\eta\,\text{ and }\eta\geq 0\Longrightarrow R[\eta]\geq 0.
$$
If $\eta\in C^2(\prt\Gw)$ satisfies $0\leq \eta\leq 1$, $\eta=1$ in a neighborhood of $K$ we take for test function 
$\gz=(R[\eta])^{p'}\phi$. There holds
$$\BA {lll}\Gd\gz=-\gl\gz+p'\phi (R[\eta])^{p'-1}\Gd R[\eta]+p'(p'-1)\phi (R[\eta])^{p'-2}|\nabla R[\eta]|^2\\[2mm]
\phantom{--------------}+2(p'-1)(R[\eta])^{p'-1}\nabla\phi.\nabla R[\eta].
\EA$$
As in \cite[Lemma 1.1]{MV01} we have
$$\BA{lll}
-\myint{\Gw}{}u\Gd\gz dx\leq \left (\myint{\Gw}{}u^p\gz dx\right)^{\frac 1p}
\left (\gl\left (\myint{\Gw}{}\gz dx\right)^{\frac {1}{p'}}+p'\left (\myint{\Gw}{}|L(\eta)|^{p'} dx\right)^{\frac {1}{p'}}\right),
\EA$$
where 
$$L(\eta)=|\phi^{\frac{1}{p'}}\Gd R[\eta]|+2|\phi^{-\frac{1}{p}}\nabla\phi.\nabla R[\eta]|.
$$
From $(\ref{D1-1})$ we have (see \cite[formula (1.2)]{MV01})
\begin{equation}\label{D6-4}\BA{lll}\left(\myint{\prt\Gw}{}\eta d\gm\right)^{p'}
+C_\gm\myint{\Gw}{}u^p\gz dx\leq mC_\gm\myint{\Gw}{}|\nabla u]^q\gz dx\\[4mm]
\phantom{--------}
+C_\gm\left (\myint{\Gw}{}u^p\gz dx\right)^{\frac 1p}
\left (\gl\left (\myint{\Gw}{}\gz dx\right)^{\frac {1}{p'}}+p'\left (\myint{\Gw}{}|L(\eta)|^{p'} dx\right)^{\frac {1}{p'}}\right),
\EA\ee
where 
$$C_\gm=\left(\myint{\prt\Gw}{}\left|\myfrac{\prt\phi}{\prt{\bf n}}\right|^{-\frac{p'}{p}}d\gm\right)^{\frac{p'}{p}}.
$$
The "optimal lifting" introduced in \cite{MV01} has the property that the mapping $\eta\mapsto L(\eta)$ is continuous from $W^{\frac{2}{p},p'}(\prt\Gw)$ into $L^{p'}(\Gw)$. Note that with $R[\eta]=\BBP_\Gw[\eta]$, which is a positive lifting, the continuity of the mapping $L$ holds only when $1<p'<2$. This is why the construction in \cite{MV01} is much more elaborate. We conclude as in Step 1 by considering a sequence $\{\eta_n\}\subset C^2(\prt\Gw)$ such that 
$0\leq \eta_n\leq 1$, $\eta_n=1$ in a neighborhood of $K$, such that $\norm{\eta_n}_{W^{\frac{2}{p},p'}}\to 0$. Then 
$\eta_n\to 0$ in $L^1(\prt\Gw)$, $\gz_n\to 0$ a.e. and $L(\eta_n)\to 0$ in $L^{p'}(\Gw)$. Thus the right-hand side of 
$(\ref{D6-4})$ tends to $0$. This ends the proof.
\qeda\medskip

\nind\Remark We conjecture that $(\ref{D6-2})$ could be strengthened and replaced by: There exists a constant $c>0$ such that 
\begin{equation}\label{D6-5}
\BA{lll}
\gm(K)\leq cCap^{\prt\Gw}_{\frac{2}{p},p'}(K)\quad\text{for any compact set } K\subset\prt\Gw.
\EA\ee
This is a necessary condition when $m=0$ (see \cite{BVHNV}).\medskip
\mysection{The boundary trace}
\subsection{The regular boundary trace}

\nind{\it Proof of \rth{Tra0}}.  Set $u=v^b$ for some $b>1$, then we have that
  \begin{equation}\label{B3+}
-\Gd v=(b-1)\myfrac{|\nabla v|^2}{v}+\myfrac{1}{b}v^{1+b(p-1)}-m b^{q-1}v^{(q-1)(b-1)}|\nabla v|^q:=F.
\ee
By H\"older's inequality,
  \begin{equation}\label{B3++}m b^{q-1}v^{(q-1)(b-1)}|\nabla v|^q\leq \myfrac{b-1}{2}\myfrac{|\nabla v|^2}{v}+m b^{q-1}\left(\myfrac{2m b^{q-1}}{b-1}\right)^{\frac{q}{2-q}}v^{\frac{2b(q-1)}{2-q}+1}.
\ee
{\it Case 1: $q< \frac{2p}{p+1}$}. There holds $\frac{2b(q-1)}{2-q}+1<1+b(p-1)$ independently of $b$. Hence for any $\gd>0$ there exists $C=C(\gd,b,m ,p,q)>0$ such that 
  \begin{equation}\label{B3+++}m b^{q-1}\left(\myfrac{2m b^{q-1}}{b-1}\right)^{\frac{q}{2-q}}v^{\frac{2b(q-1)}{2-q}+1}\leq \frac \gd b v^{1+b(p-1)}+C.
\ee
Therefore 
  \begin{equation}\label{B4}\BA{lll}
F\geq \myfrac{b-1}{2}\myfrac{|\nabla v|^2}{v}+\myfrac{1-\gd}{b}v^{1+b(p-1)}-C.
\EA\ee
If $\psi=\BBG_{B_R}[1]$ (ie. the solution of $-\Gd\psi=1$ in $B_R$ vanishing on $\prt B_R$), we have
$$-\Gd (v+C\psi)\geq \myfrac{b-1}{2}\myfrac{|\nabla v|^2}{v}+\myfrac{1-\gd}{b}v^{1+b(p-1)}\geq 0.
$$By Doob's theorem  on positive superharmonic functions (see \cite{Doo}) we have that $\myfrac{|\nabla v|^2}{v}+v^{1+b(p-1)}\in L^1_\gr(\Gw)$. We put $a=b^{-1}-1$, then  $a<0$ and 
$v=u^\frac 1b=u^{1+a}$. Therefore 
$$\nabla v=(1+a)u^a\nabla u\,,\;\frac{|\nabla v|^2}{v}=(1+a)^2u^{a-1}|\nabla u|^2\text{ and }
v^{1+b(p-1)}=u^{p+a},
$$
consequently
$$u^{a-1}|\nabla u|^2+u^{p+a}\in L^1_\gr(\Gw).
$$
Let $1<\ell<\frac{2p}{p+1}<2$, then 
\bel{D03}\BA{lll}
\myint{\Gw}{}|\nabla u|^\ell\gr dx=\myint{\Gw}{}|u^{\frac{a-1}{2}}\nabla u|^\ell u^{\frac{(1-a)\gl}{2}}\gr dx\\[4mm]
\phantom{\myint{\Gw}{}|\nabla u|^\ell\gr dx}
\leq \ge\myint{\Gw}{}u^{a-1}|\nabla u|^2\gr dx+C(\ge)\myint{\Gw}{}u^{\frac{(1-a)\ell}{2-\ell}}\gr dx.
\EA\ee
We fix $a<0$ such that $\frac{(1-a)\ell}{2-\ell}=p+a$, or equivalently 
$$a=-\frac{p+1}{2}\left(\frac{2p}{p+1}-\ell\right).
$$
Finally, we infer that for any $\ell<\frac{2p}{p+1}$, $|\nabla u|^\ell\in L^1_\gr(\Gw)$. This implies in particular that 
$|\nabla u|^q\in L^1_\gr(\Gw)$. \smallskip

Let $\Psi=m\BBG_\Gw[|\nabla u|^q]$, then $\Psi>0$ and
$$-\Gd(u+\Psi)=u^p.
$$
Clearly the function $u+\Psi$ is positive and superharmonic in $\Gw$. By using again Doob's theorem  \cite{Doo}, it follows that 
$[-\Gd(u+\Psi)]=u^p\in L^1_\gr(\Gw)$ and there exists a nonnegative Radon measure $\gm$ on $\prt\Gw$ such that 
$$u=\BBG_\Gw[u^p]-\Psi+\BBP_\Gw[\gm]=\BBG_\Gw[u^p-m|\nabla u|^q|]+\BBP_\Gw[\gm],
$$
where $\BBP_\Gw$ is the Poisson operator in $\Gw$. This implies that $(\ref{D01})$ holds.\medskip

\nind{\it Case 2: $q=\frac{2p}{p+1}$}. We proceed as in the proof of \rth{Tra0}, setting $u=v^b$, $b>1$. Since 
$q=\frac{2p}{p+1}$, inequality $(\ref{B3++})$ becomes
\begin{equation}\label{B8+}m b^{\frac{p-1}{p+1}}v^{\frac{(p-1)(b-1)}{p+1}}|\nabla v|^{\frac{2p}{p+1}}\leq \myfrac{b-1}{2}\myfrac{|\nabla v|^2}{v}+
m b^{\frac{p-1}{p+1}}\left(\myfrac{2m b^{\frac{p-1}{p+1}}}{b-1}\right)^{p}v^{1+b(p-1)}.
\ee
Defining $m_1$ by the identity, 
\bel{B9}
m_1=\left(\myfrac{b-1}{2b}\right)^{\frac p{p+1}},
\ee
we deduce that for $0<m<m_1$ and some $\gd\in (0,1)$, there holds
\bel{D04}
-\Gd v\geq \myfrac{b-1}{2}\myfrac{|\nabla v|^2}{v}+\myfrac{b-1}{2}v^{1+b(p-1)}.
\ee
Again, by  Doob's theorem, $\frac{|\nabla v|^2}{v}+v^{1+b(p-1)}\in L^1_\gr(\Gw)$, which implies that $\sqrt v\in W^{1,2}_\gr(\Gw)$. Using Sobolev type imbedding theorem for weighted Sobolev spaces (see e.g. \cite[Section 19]{OpKu}), 
\bel{D05}
\left(\myint{\Gw}{}(\sqrt v)^{\frac{2(N+1)}{N-1}}\gr dx\right)^{\frac{N-1}{N_1}}\leq c
\myint{\Gw}{}\left((\sqrt v)^{2}+|\nabla\sqrt v|^2\right)\gr dx.
\ee
If we choose in particular $b=\frac{N+1}{p(N-1)}$ we deduce that $u^p\in L^1_\gr(\Gw)$. Actually, for any $1\leq \tilde p<\frac{N+1}{N-1}$,
$u^{\tilde p}\in L^1_\gr(\Gw)$ and for any $\ge>0$, $\frac{|\nabla u|^2}{u^{1+\ge}}\in L^1_\gr(\Gw)$. 
We have from $(\ref{D03})$ with $\ell=\frac{2p}{p+1}$, 
\bel{D06}\BA{lll}
\myint{\Gw}{}|\nabla u|^\frac{2p}{p+1}\gr dx=\myint{\Gw}{}u^{-\frac{(1+\ge)p}{p+1}}|\nabla u|^\frac{2p}{p+1} u^{\frac{(1+\ge)p}{p+1}}\gr dx\\[4mm]
\phantom{\myint{\Gw}{}|\nabla u|^\frac{2p}{p+1}\gr dx}
\leq s\myint{\Gw}{}\myfrac{|\nabla u|^2}{u^{1+\ge}}\gr dx+C(s)\myint{\Gw}{}u^{(1+\ge)p}\gr dx.
\EA\ee
If $\ge$ is chosen such that $(1+\ge)p=\tilde p<\frac{N+1}{N-1}$, we infer that $|\nabla u|^\frac{2p}{p+1}\in L^1_\gr(\Gw)$. We end the proof as in Case 1.\qeda\medskip

\nind\Remark The same regularity and boundary trace results hold if it is assumed that
$u$ is a nonnegative supersolution of $\ref{An1}$ in $\Gw$.

\subsection{The singular boundary trace}

{\it Proof of \rth{obs1}}.  {\it Assertion 1}. We assume that $F:=|\nabla u|^q\in L^1_\gr(B_\ge(z)\cap\Gw)$. We set $ F_\ge=F{\bf 1}_{_{B_\ge(z)\cap\Gw}}$ and 
$\Psi_\ge=\BBG_{B_\ge(z)\cap\Gw}[F_\ge]$. Then $\Psi_\ge$ has boundary trace zero on $B_\ge\cap\prt\Gw$ and
$$-\Gd(u+m\Psi_\ge)=u^p\quad\text{in }B_\ge(z)\cap\Gw.$$
Thus $u+m\Psi_\ge$ is a positive super-harmonic function in $B_\ge\cap\Gw$. Hence $u^p\in L^1_\gr(B_\ge(z)\cap\Gw)$ and there exists a Radon measure $\gm_\ge$ such that 
$u+m\Psi_\ge$ admits for boundary trace $\gm_\ge$ on $B_\ge(z)\cap\prt\Gw$. This implies that $u$ admits the same boundary trace on $B_\ge(z)\cap\prt\Gw$. \smallskip
{\bf 1}
\nind {\it Assertion 2}. We assume that $H:=u^p{\bf 1}_{_{B_\ge(z)\cap\Gw}}\in L^1_\gr(B_\ge(z)\cap\Gw)$. If $F_\ge=|\nabla u|^q{\bf 1}_{_{B_\ge(z)\cap\Gw}}\in L^1_\gr(B_\ge(z)\cap\Gw)$, we deduce from Assertion 1 that $u$ admits the boundary trace $\gm_\ge\in\mathfrak M_+(B_\ge(z)\cap\prt\Gw)$ on $B_\ge(z)\cap\prt\Gw$. 
If for any $\ge'\in (0,\ge]$
$$\myint{B_{\ge'}(z)\cap\Gw}{}|\nabla u|^q\gr dx=\infty,
$$
there holds
$$\myint{B_{\ge'}(z)\cap\Gw}{}\left(m|\nabla u|^q-u^p\right)\gr dx=\infty.
$$
For $0<\gd<\frac{\ge'}{2}$, set $\Gth_{\gd,\ge'}=B_{\ge'}(z)\cap\Gw\cap\{x\in\Gw:\gr(x)>\gd\}$ and denote by $\phi_{\gd,\ge'}$ the first eigenfunction of $-\Gd$ in $H^1_0(\Gth_{\gd,\ge'})$ normalized by $\sup \phi_{\gd,\ge'}=1$ and let $\gl_{\gd,\ge'}$ be the corresponding eigenvalue. Then $\phi_{\gd,\ge'}\to \phi_{0,\ge'}$, uniformly,  $\gl_{\gd,\ge'}\downarrow\gl_{0,\ge'}$ and $\frac{\prt \phi_{\gd,\ge'}}{\prt{\bf n}}\to \frac{\prt \phi_{0,\ge'}}{\prt{\bf n}}$ in the sense that 
$$\frac{\prt \phi_{\gd,\ge'}}{\prt{\bf n}}(x+\gd{\bf n})\to \frac{\prt \phi_{0,\ge'}}{\prt{\bf n}}(x)\quad\text{uniformly for }x\in \prt\Gw\cap  B_{\ge'}(z).$$
Let $v_{\ge',\gd}$ be the solution of 
$$\BA{lll}
-\Gd v+m|\nabla u|^q-u^p=0\quad&\text{in }\Gth_{\gd,\ge'}\\[1mm]
\phantom{-\Gd +m|\nabla u|^q-u^p}v=u&\text{on }\prt\Gth^{up}_{\gd,\ge'}:=\overline\Gth_{\gd,\ge'}\cap\{x:\gr(x)=\gd\}\\[1mm]
\phantom{-\Gd +m|\nabla u|^q-u^p}v=0&\text{on }\prt\Gth^{lat}_{\gd,\ge'}:=\prt\Gth_{\gd,\ge'}\cap\{x:\gr(x)>\gd\}.
\EA$$
Then $u\geq v_{\ge',\gd}$ in $\overline\Gth_{\gd,\ge'}$ and 
\bel{D07}
\myint{\Gth_{\gd,\ge'}}{}\left(\gl_{\gd,\ge'}v+m|\nabla u|^q-u^p\right)\phi_{\gd,\ge'} dx=-\myint{\prt\Gth^{up}_{\gd,\ge'}}{}
\frac{\prt \phi_{\gd,\ge'}}{\prt{\bf n}} udS.
\ee
Since the left-hand side of $(\ref{D07})$ tends to $\infty$ when $\gd\to 0$, we deduce that 
\bel{D08}\displaystyle
\lim_{\gd\to 0}\myint{\Gw\cap B_{\ge'}(z)}{}
udS=\infty.
\ee
Thus $z\in \CS(u)$.
\qeda\medskip

\nind\Remark Note also that if $p>2$, then $u^p\in L^1_\gr(\Gw)$ implies $u\in L^1(\Gw)$ and the assertion 2 follows from \cite[Lemma 2.8]{MV001}. If $p>\frac{N+1}{N-1}$ and if we assume that $u$ satisfies 
\begin{equation}\label{D11}
u(x)\leq c(\gr(x))^{-\frac{2}{p-1}},
\ee
then $u^p\in L^1_\gr(\Gw)$. 
\medskip

In order to describe the boundary singularities of solutions we introduce the following equation studied in \cite{NGT-LV}

\begin{equation}\label{D12}
\BA{lll}\displaystyle
-\Gd'\chi-\gb(\gb+2-N)\chi+m\left(\gb^2\chi^2+|\nabla'\chi|^2\right)^{\frac q2}=0\quad&\text{in } S_+^{N-1}\\
\phantom{-\Gd'-\gb(\gb+2-N)\chi+m\left(\gb^2\chi^2+|\nabla'\chi|^2\right)^{\frac q2}}
\chi=0\quad&\text{in } \prt S_+^{N-1},
\EA\ee
 where $m>0$ and 
 \begin{equation}\label{D12-0}
 \gb=\frac{2-q}{q-1}
 \ee
It is proved in \cite{NGT-LV} that if $1<q<\frac{N+1}{N}$, $(\ref{D12})$ admits a unique solution $\chi_{_m}$. The function $V_{\chi_{_m}}(x)=V_{\chi_{_m}}(r,s)=r^{-\gb}\chi_{_m}(s)$ where $(r,s)\in \BBR_+\ti S_+^{N-1}$ is the only positive solution of  
 \begin{equation}\label{D13}
 \BA{lll}\displaystyle
-\Gd v+m|\nabla v|^q=0\quad&\text{in } \BBR_+^{N}
\EA\ee
 which vanishes on $\prt\BBR_+^{N}\setminus\{0\}$ and satisfies
   \begin{equation}\label{D14}
 \BA{lll}\displaystyle\limsup_{x\to 0}|x|^{N-1} v(x)=\infty.
\EA\ee
It is a consequence of uniqueness that 
   \begin{equation}\label{D14m}
 \BA{lll}\chi_{_m}=m^{-\frac{1}{q-1}}\chi_1:=m^{-\frac{1}{q-1}}\chi.
\EA\ee
Furthermore, if $v_{k\gd_0}$ is the unique positive solution of 
 \begin{equation}\label{D15}
 \BA{lll}\displaystyle
-\Gd v+m|\nabla v|^q=0\quad&\text{in } \BBR_+^{N}\\
\phantom{-\Gd +m|\nabla v|^q}
v=k\gd_0&\text{on } \prt\BBR_+^{N},
\EA\ee
then $v_{k\gd_0}\uparrow v_{\chi_{_m}}$ when $k\to\infty$. If $\BBR^N_+$ is replaced by a bounded smooth subset $\Gw$, the previous statements still hold provided some adaptations are performed. We assume that $\Gw$ is in normal position with respect to $0\in\prt\Gw$. The next result is proved in \cite{NGT-LV}.
\bth{NgVTh} Let $\Gw$ be as described above, $m>0$ and $1<q<\frac{N+1}{N}$. \\
1- Then for any $k>0$ there exists a unique positive weak solution $v_{k\gd_0}$ of 
 \begin{equation}\label{D16}
 \BA{lll}\displaystyle
-\Gd v+m|\nabla v|^q=0\quad&\text{in } \Gw\\
\phantom{-\Gd +m|\nabla v|^q}
v=k\gd_0&\text{on } \prt\Gw.
\EA\ee
Furthermore 
 \begin{equation}\label{D17}
\lim_{\tiny\BA{lll}\Gw\ni x\to 0\\
\frac{x}{|x|}\to s\EA}|x|^{N-1}v_k(x)=c_Nk\phi_1( s)\quad\text{locally uniformly in } s\in S^{N-1}_+.
\ee
2- The function $v_k$ is stable in the sense that if $\{\gm_n\}$ is a sequence of positive Radon measures on $\prt\Gw$ which converges weakly to $k\gd_0$, then the corresponding sequence of solutions $\{v_{\gm_n}\}$ of 
 \begin{equation}\label{D18}
 \BA{lll}\displaystyle
-\Gd v+m|\nabla v|^q=0\quad&\text{in } \Gw\\
\phantom{-\Gd +m|\nabla v|^q}
v=\gm_n&\text{on } \prt\Gw,
\EA\ee
converges locally uniformly in $\Gw$ to  $v_{k\gd_0}$. \\
3- Finally, when $k\uparrow\infty$, $v_{k\gd_0}\uparrow v_{\chi_{_m}}$ where $v_{\chi_{_m}}$ is the unique positive solution of 
 \begin{equation}\label{D20}
 \BA{lll}\displaystyle
-\Gd v+m|\nabla v|^q=0\quad&\text{in } \Gw,
\EA\ee
which vanishes on $\prt\Gw\setminus\{0\}$ and satisfies $(\ref{D14})$. Furthermore  $v_{\chi_{_m}}$ verifies the following limits, locally uniformly on 
$S^{N-1}_+$,
 \begin{equation}\label{D21}
\lim_{\tiny\BA{lll}x\in\Gw\\x\to 0\\
\frac{x}{|x|}\to s\EA}
|x|^{\gb}v_{\chi_{_m}}(x)={\chi_{_m}}( s),
\ee
and
 \begin{equation}\label{D21-0}\BA{lll}
\displaystyle\lim_{\tiny\BA{lll}x\in\Gw\\x\to 0\\
\frac{x}{|x|}\to s\EA}
|x|^{\gb+1}\frac{x}{|x|}.\nabla v_{\chi_{_m}}(x)&\!\!\!\!=-\gb{\chi_{_m}}( s),\\
\phantom{}
\displaystyle\lim_{\tiny\BA{lll}x\in\Gw\\x\to 0\\
\frac{x}{|x|}\to s\EA}
|x|^{\gb+1}\nabla_{tang} v_{\chi_{_m}}(x)&\!\!\!\!=\nabla_{tang}{\chi_{_m}}( s),
\EA\ee
where $\nabla_{tang}=r^{-1}\nabla '$ denotes the tangential gradient.
\es

\nind We set 
 \begin{equation}\label{D22}
\Gw'_\gd=\{x\in\Gw:\gr(x)>\gd\}\,,\; \Gw_\gd=\{x\in\Gw:0<\gr(x)<\gd\}\,\text{ and }\;\Gs_\gd=\prt\Gw'_\gd.
\ee
It is known that $\Gs_\gd$ is smooth for $\gd$ small enough.  The following variant is proved in \cite[Corollary 2.4]{NGT-LV}. 
\bcor{NgVTh} Under the assumptions on $N$, $q$ and $m$ of \rth{NgVTh}, assume that $\{\gd_n\}$ is a sequence decreasing to $0$, $\{\gm_n\}$ is a sequence of positive bounded Radon measures on $\Gs_{\gd_n}$ which converges in the sense of measures in $\overline\Gw$ to a measure $\gm$ on $\prt\Gw$. Then the sequence $\{v_{\gm_n}\}$ of solutions of 
 \begin{equation}\label{D23}
 \BA{lll}\displaystyle
-\Gd v+m|\nabla v|^q=0\quad&\text{in } \Gw'_{\gd_n}\\
\phantom{-\Gd +m|\nabla v|^q}
v=\gm_n&\text{on } \Gs_{\gd_n},
\EA\ee
converges up to a subsequence locally uniformly in $\Gw$ to a positive solution $v_{\gm}$ of 
 \begin{equation}\label{D24}
 \BA{lll}\displaystyle
-\Gd v+m|\nabla v|^q=0\quad&\text{in } \Gw\\
\phantom{-\Gd +m|\nabla v|^q}
v=\gm &\text{on } \prt\Gw.
\EA\ee
\es

\bprop{obsL2} Let $p>1$, $1<q<\frac{N+1}{N}$ and $m>0$. Let $u$ be a positive solution of  $(\ref{An1})$ in $\Gw$ such that there exist a sequence $\{z_n\}\subset\prt\Gw$ converging to $z$ and two decreasing sequences $\{\ge_n\}$ and $\{\gd_n\}$ converging to $0$ such that 
\begin{equation}\label{D25}
\BA{lll}\displaystyle
\lim_{n\to\to \infty}\myint{B_{\ge_n}(z_n)\cap\Gs_{\gd_n}}{}udx=\infty, 
\EA\ee
then there holds
\begin{equation}\label{D26}
\BA{lll}\displaystyle
\liminf_{\tiny\BA{lll}x\in\Gw\\ x\to z\\
\!\!\frac{x-z}{|x-z|}\to s\EA}|x-z|^{\frac{2-q}{q-1}}u(x)\geq\chi ( s)\quad\text{locally uniformly in } s\in S^{N-1}_+.
\EA\ee
\es
\Proof For $k>0$, there exists $n_0$ such that for $n\geq n_0$, 
$$\myint{B_{\ge_n}(z_n)\cap\Gs_{\gd_n}}{}udx>k. 
$$
 Hence there exists $\ell:=\ell_n>0$ such that
 $$\myint{B_{\ge_n}(z_n)\cap\Gs_{\gd_n}}{}\min\{u,\ell\}dx=k. 
$$
We set $\gm_{n,k}=\min\{u,\ell\}\lfloor_{\Gs_{\gd_n}}{\bf 1}_{_{\Gs_{\gd_n}\cap B_{\ge_n}(z_n)}}$ and denote by $v_{\gm_{n,k}}$ the corresponding solution of $(\ref{D23})$ in $\Gw'_{\gd_n}$. Then $u\geq v_{\gm_{n,k}}$ in $\Gw'_{\gd_n}$. Up to a rotation we can assume that $\prt\BBR^N_+$ is tangent to $\prt\Gw$ at $z$. Using \rcor{NgVTh} we obtain $u\geq v_{k\gd_z}$. Letting $k\to\infty$ and using \rth{NgVTh}-3 we deduce that 
\begin{equation}\label{D27}
\liminf_{\tiny\BA{lll}x\in\Gw\\ x\to z\\
\!\!\frac{x-z}{|x-z|}\to s\EA}|x-z|^{\gb}u(x)\geq \chi( s)\quad\text{locally uniformly in } s\in S^{N-1}_+.
\ee
\qeda \medskip

In the sequel we denote $\chi_1=\chi$ and $v_{\chi_1}=v_{\chi}$. In the next theorem we show the existence of positive singular solution of ($\ref{An4}$) with a strong blow-up in $|x|^{-\gb}$ provided 
the function $v_{\chi_1}$ has no critical point in $\Gw$ and $\frac{2p}{p+1}q<\frac{N+1}{N}$.  If it is the case the constant $M_{v_{\chi}}$ defined below is positive because of $(\ref{D21-0})$ and Hopf boundary lemma,
\begin{equation}\label{D28}\displaystyle
M_{v_{\chi}}=\min_{x\in\Gw}\myfrac{|\nabla v_{\chi}(x)|^q}{v_{\chi}^p(x)}.
\ee

\bth{singbeta} Let $\Gw$ be a bounded smooth domain with $0\in\prt\Gw$,  $p>1$ and $\frac{2p}{p+1}<q<\frac{N+1}{N}$. If 
\begin{equation}\label{D29}
m>m_{v_{\chi}}=:\myfrac{p-1}{p-q}\left(\myfrac{p-q}{(q-1)M_{v_{\chi}}}\right)^{\frac{q-1}{p-1}},
\ee
then there exists a positive solution of ($\ref{An4}$) which satisfies 
\begin{equation}\label{D30}
\lim_{\tiny\BA{lll}x\in\Gw\\ x\to 0\\
\!\!\frac{x-z}{|x-z|}\to s\EA}|x|^{\gb}u(x)= \chi_{_m}( s)\quad\text{locally uniformly in } s\in S^{N-1}_+,
\ee
where $\chi_{_m}$ is the unique positive solution of ($\ref{D12}$).
\es
\Proof The function $\chi$ is the unique positive solution of ($\ref{D12}$), and since it depends on $m>0$, we denote it by $\chi_{_m}$. Clearly $\chi_{_m}=m^{-\frac{1}{q-1}}\chi$. Then $v_{\chi_{_m}}=m^{-\frac{1}{q-1}}v_{\chi}$ is the solution of ($\ref{D20}$) which is obtained in \rth{NgVTh}, since this solution is the unique positive solution of
($\ref{D20}$) which satisfies ($\ref{D21}$)-($\ref{D21-0}$). 
We also set  
$$L_{m,p,q}u=-\Gd u+m|\nabla u|^q-u^p. 
$$
The function $v_{\chi_{_m}}$ is a subsolution of ($\ref{An4}$).  Let $0<\tilde m< m$, then $v_{\chi_{_m}}<v_{\chi_{\tilde m}}$. Furthermore 
$$\BA{lll}L_{m,p,q}v_{\chi_{\tilde m}}=(m-\tilde m)|\nabla v_{\chi_{\tilde m}}|^q- v^p_{\chi_{\tilde m}}\\[2mm]\phantom{L_{m,p,q}v_{\chi_{\tilde m}}}
=(m-\tilde m)\tilde m^{-\frac{q}{q-1}}|\nabla v_{\chi}|^q-\tilde m^{-\frac{p}{q-1}}v^p_{\chi}\\[2mm]\phantom{L_{m,p,q}v_{\chi_{\tilde m}}}
\geq \left((m-\tilde m)\tilde m^{-\frac{q}{q-1}}M_{v_{\chi}}-\tilde m^{-\frac{p}{q-1}}\right)v^p_{\chi}
\\[2mm]\phantom{L_{m,p,q}v_{\chi_{\tilde m}}}
\geq \left(m-\left(\tilde m+\myfrac{1}{\tilde m^{\frac{p-q}{q-1}}M_{v_{\chi}}}\right)\right)\tilde m^{-\frac{q}{q-1}}M_{v_{\chi}}v^p_{\chi}.
\EA$$
Then
\begin{equation}\label{D31}
\displaystyle\min_{X>0}\left\{X+\myfrac{1}{X^{\frac{p-q}{q-1}}M_{v_{\chi}}}\right\}=\myfrac{p-1}{p-q}\left(\myfrac{p-q}{(q-1)M_{v_{\chi}}}\right)^{\frac{q-1}{p-1}}:=m_{v_{\chi}}
\ee
and the minimum is achieved for
\begin{equation}\label{D32}
X=X_0=\left(\myfrac{p-q}{(q-1)M_{v_{\chi}}}\right)^{\frac{q-1}{p-1}}.
\ee
If we fix $\tilde m=X_0$ it follows that for $m>m_{v_{\chi}}$, the function $v_{\chi_{\tilde m}}$ satisfies $L_{m,p,q}v_{\chi_{\tilde m}}\geq 0$ in $\Gw$ and it is larger than the subsolution $v_{\chi_{ m}}$. Hence there exists a solution $u$ of ($\ref{An4}$) in $\Gw$ and it satisfies 
\begin{equation}\label{D33}
v_{\chi_{ m}}\leq u\leq v_{\chi_{\tilde m}}\quad\text{in }
\Gw.\ee
The end of the proof is standard. For $\ell>0$ we set $S_\ell[v](x)=\ell^\gb v(\ell x)$. Then $u_\ell:=S_\ell[u]$ satisfies 
\begin{equation}\label{D34}
-\Gd u_\ell+m|\nabla u_\ell|^q-\ell^\frac{q(p+1)-2p}{q-1}u_\ell^p=0\quad\text{in }
\Gw_\ell:=\tfrac 1\ell\Gw,\ee
and 
$$S_\ell[v_{\chi_{ m}}]\leq u_\ell\leq S_\ell[v_{\chi_{\tilde m}}]\quad\text{in }
\Gw_\ell.
$$
By \rth{BS1}, 
\begin{equation}\label{D35}
|\nabla u_\ell(x)|+\frac{\ell u_\ell(x)}{\gr(\ell x)}\leq c|x|^{-\gb-1}\quad\text{in }\overline \Gw_\ell\setminus\{0\}.
\ee
Since $\prt\Gw$ is smooth, there exists $\ge_0>0$ such that $c_2\ell\gr_\ell(x)\leq \gr(\ell x)\leq c_1\ell\gr_\ell(x)$ for $|x|\leq\ge_0$, in which formula
we denote $\gr_\ell(x)=\dist(x,\Gw_\ell)$. Since $q(p+1)-2p>0$, $\ell^\frac{q(p+1)-2p}{q-1}u_\ell^p\to 0$ when $\ell\to 0$, locally uniformly in $\Gw_\ell\cap B_\gd^c$ for any $\gd>0$ and by standard elliptic equations regularity results \cite{GT},  $D^2u_\ell$ is also locally bounded in $\Gw_\ell\cap B_\gd^c$. When $\ell\to 0$, $S_\ell[v_{\chi_{ m}}]$ and $S_\ell[v_{\chi_{\tilde m}}]$ converge respectively to $x\mapsto |x|^{-\gb}\chi_{_m}(\frac{x}{|x|})$ and $x\mapsto |x|^{-\gb}\chi_{\tilde m}(\frac{x}{|x|})$. Therefore, if $\displaystyle u=\lim_{n\to\infty}u_{\ell_n}$ for some sequence $\{\ell_n\}$ converging to $0$, the function $u$ is nonnegative and satisfies 
\begin{equation}\label{D36}
-\Gd u+m|\nabla u|^q=0
\ee
in  $\BBR^N_+$ and there holds
$$|x|^{-\gb}\chi_{_m}(\frac{x}{|x|})\leq u(x)\leq |x|^{-\gb}\chi_{\tilde m}(\frac{x}{|x|}).
$$
Since $(\ref{D36})$ admits a unique positive solution vanishing on $\prt\BBR^N_+\setminus\{0\}$ such that $\displaystyle \limsup_{x\to 0}|x|^\gb u(x)>0$ (see \cite[Proposition 3.24-Step 2 ]{NGT-LV}), it follows that $u(x)=|x|^{-\gb}\chi_{_m}(\frac{x}{|x|})$. Uniqueness implies that $u_{\ell}\to u$ and $(\ref{D30})$ holds.\qeda\medskip

\nind \Remark The assumption that $v_{\chi}$ admits no critical point in $\Gw$ is uneasy to verify. At least it is easy to see that $v_{\chi}$ cannot have any non-degenerate critical point in $\Gw$. Furthermore, because of Hopf boundary lemma and the behaviour of  $v_{\chi}$ near $x=0$ given by $(\ref{D21})$, the critical points of  $v_{\chi}$ are located in a compact subset $N$ of $\Gw$, possibly empty. For $\ge>0$ we set 
$$N_\ge=\{x\in\Gw:\dist(x,N)<\ge\}.
$$
 If $\ge$ is small enough $\overline N_\ge\subset\Gw$. Denote 
 \begin{equation}\label{D37}\displaystyle
M^\ge_{v_{\chi}}=\min_{x\in\Gw\setminus N_\ge}\myfrac{|\nabla v_{\chi}(x)|^q}{v_{\chi}^p(x)}\quad{and }\; \;
m^\ge_{v_{\chi}}=:\myfrac{p-1}{p-q}\left(\myfrac{p-q}{(q-1)M^\ge_{v_{\chi}}}\right)^{\frac{q-1}{p-1}}.
\ee
The proof of the next result is similar to the one of \rth{singbeta}.

\bth{singbeta2} Let $\Gw$ be a bounded smooth domain with $0\in\prt\Gw$,  $p>1$ and $\frac{2p}{p+1}<q<\frac{N+1}{N}$. If 
$N$ denotes the set of critical points of $v_{\chi}$, then for any $\ge>0$ small enough and $m>m^\ge_{v_{\chi}}$
there exists a positive solution of ($\ref{An1}$) in $\Gw\setminus N_\ge$ which vanishes on $\prt\Gw$ and satisfies $(\ref{D30})$.
\es

\end{document}